\documentclass[11pt]{article}
\usepackage{amsmath,amssymb,latexsym,float}
\usepackage{latexsym}
\usepackage{amsthm}
\usepackage{cases}

\usepackage{scalefnt}
\usepackage{booktabs}
\usepackage[english]{babel}
\usepackage{graphicx}

\usepackage{titlesec}
\usepackage{epsfig}
\usepackage{enumerate}
\usepackage{subcaption}
\usepackage{booktabs}
\usepackage{natbib}
\usepackage{graphicx}
\usepackage{multirow}
\usepackage{rotating}
\usepackage{comment}

\usepackage{xargs}                      
\usepackage[pdftex,dvipsnames]{xcolor}  
\usepackage[colorinlistoftodos,prependcaption,textsize=small]{todonotes}
\newcommandx{\unsure}[2][1=]{\todo[linecolor=red,backgroundcolor=red!25,bordercolor=red,#1]{#2}}
\newcommandx{\change}[2][1=]{\todo[linecolor=blue,backgroundcolor=blue!25,bordercolor=blue,#1]{#2}}
\newcommandx{\info}[2][1=]{\todo[linecolor=OliveGreen,backgroundcolor=OliveGreen!2,bordercolor=OliveGreen,#1]{#2}}
\newcommandx{\improvement}[2][1=]{\todo[linecolor=Plum,backgroundcolor=Plum!25,bordercolor=Plum,#1]{#2}}
\newcommandx{\thiswillnotshow}[2][1=]{\todo[disable,#1]{#2}}

\usepackage{fancyhdr}
\usepackage{setspace}
\usepackage[section]{placeins}

\usepackage{xcolor}
\usepackage{amsmath,algorithm,algpseudocode}
\makeatletter
\renewcommand{\Function}[2]{%
  \csname ALG@cmd@\ALG@L @Function\endcsname{#1}{#2}%
  \def\jayden@currentfunction{#1}%
}
\newcommand{\funclabel}[1]{%
  \@bsphack
  \protected@write\@auxout{}{%
    \string\newlabel{#1}{{\jayden@currentfunction}{\thepage}}%
  }%
  \@esphack
}
\makeatother



\usepackage{breqn}
\usepackage{gensymb}
\usepackage{array, makecell}
\usepackage{pseudocode}
\usepackage{algorithmicx}
\usepackage{algorithm}
\usepackage{url} 
\usepackage{appendix}
\usepackage{ragged2e}
\usepackage{longtable}
\usepackage{multirow}
\usepackage{physics}
\usepackage[]{natbib} 
\newtheorem{proposition}{Proposition}
\newtheorem{definition}{Definition}

\usepackage[nocomma]{optidef}
\usepackage{mathtools,zref-savepos}
\newtagform{roman}[]()

\topmargin by -\headsep \textheight 8.9in \oddsidemargin 0pt
\evensidemargin \oddsidemargin \marginparwidth 0.5in \textwidth
6.5in
\date{}
\usepackage{tabularx}
\usepackage{graphicx}
\usepackage{bbm}
\usepackage{mathtools}
\usepackage{graphicx}
\usepackage{caption}
\usepackage{subcaption}

\usepackage{algorithm}
\usepackage{algpseudocode}
\allowdisplaybreaks
\usepackage{mathrsfs}
\usepackage{tikz}
\usetikzlibrary{decorations.pathreplacing}
\usetikzlibrary{backgrounds}
\usetikzlibrary{arrows,shapes,positioning}
\DeclareMathOperator*{\argmin}{argmin}

\usepackage{calligra}
\tikzset{el1/.style={draw, circle, fill=white, inner sep=2pt, font=\small}}
\begin{document}

\title{Repair Crew Routing for Infrastructure Network Restoration \\under Incomplete Information}
\author{}
\date{}

\author{
Subhojit Biswas \thanks{subhojit.biswas@tamu.edu} \and
Bahar \c{C}avdar \thanks{cavdab2@rpi.edu} \and Joseph Geunes \thanks{geunes@tamu.edu}}
\date{%
     $^*$$^\ddagger$\small{Department of Industrial and Systems Engineering, Texas A\&M University, College Station, TX 77843}\\[2ex]%
    $^\dagger$\small{Department of Industrial and Systems Engineering, Rensselaer Polytechnic Institute, Troy, NY, 12180}%
}


\maketitle

\begin{abstract}
{
This paper considers a disrupted infrastructure network where the repair crew knows the locations of service outages but not the locations of actual faults. Our goal is to determine a route for a single crew to visit and repair the disruptions to restore service with minimum negative impact. We call this problem the Traveling Repairman Network Restoration Problem (TRNRP). This problem presents strong computational challenges due to the combinatorial nature of the decisions, inter-dependencies within the underlying infrastructure network, and incomplete information. Considering the dynamic nature of the decisions as a result of dynamic information revelation on the status of the nodes,  we model this problem as a finite-horizon Markov decision process. Our solution approach uses value approximation based on reinforcement learning, which is strengthened by structural results that identify a set of suboptimal moves. In addition, we propose state aggregation methods to reduce the size of the state space. We perform extensive computational studies to characterize the performance of our solution methods under different parameter settings and to compare them with benchmark solution approaches. 

}
\end{abstract}

\doublespacing

\section{Introduction}

The vulnerability of critical infrastructure networks coupled with increasing extreme weather phenomena, such as hurricanes, storms and wildfires, pose substantial challenges to the functionality of these vital networks and can lead to significant negative impacts. This was starkly evident during winter storm Uri in February 2021, where approximately 4.5 million households in Texas endured multiple days of power outages (as reported by \cite{pollock2021gov} and \cite{busby2021cascading}), causing interruptions in water supply and medical services, as highlighted by \cite{mcnamara2021over}. According to AccuWeather, the economic cost due to lost productivity and damages reached an estimated \$130 billion in Texas alone and \$150 billion nationwide during the winter of 2021 \citep{puelo2021damages}. To mitigate the societal and economic consequences of extended service disruptions, timely restoration efforts are crucial to enhance communities' overall resilience in the aftermath of disruptions. 
Inadequate preparation and response decisions can impede the restoration process and significantly extend disruption periods. 
Motivated by the negative impact of service failures due to infrastructure network disruptions, in this paper, we study the logistics of network restoration decisions in the aftermath of a disaster.

We consider the aftermath of an event in which a power distribution network is disrupted, causing a loss of power over a service region. In such cases, it is not uncommon for a repair crew dispatcher to have information on the locations of power outages as observed by customers, but not necessarily the locations of the actual disruptions, since this may require inspection of the network.
This phenomenon is a result of network inter-dependencies. When a fault occurs at a power node, downstream locations in the network also lose service regardless of whether or not they experience a fault. Therefore, the repair crew often has incomplete information on the actual fault locations. 
As the duration of the service outage increases, the negative impacts on customers increase as well. Our goal is to determine a node visitation and repair sequence for a single crew that minimizes the negative impact of power loss as measured by the total service outage duration experienced by all demand locations. This problem has similarities to the Traveling Repairman Problem (TRP). However, in the TRP, a demand location's service is completed after a visit to the location. On the other hand, this is not the case in infrastructure network restoration. For a node's service to be completed, power must be restored, i.e., all faults between the node and the power source must be eliminated. This aspect of the problem increases the computational burden required to compute good solutions, posing challenges, especially in real-time decision-making in the aftermath of a disaster.

Motivated by the practical implications and the computational challenges of the routing decisions in infrastructure network restoration, in this paper, we study repair crew routing for power network restoration under incomplete information. We call the underlying problem the Traveling Repairman Network Restoration Problem (TRNRP). We formulate this problem as a Markov decision process (MDP), where the revelation of the true fault status of nodes is modeled through state transitions. To solve the problem more efficiently, we perform structural analysis to identify a set of suboptimal routes based on available road and power network information. We propose a problem-specific reinforcement learning-based approach to implement the MDP model. 
We observe that due to the inter-dependencies of the infrastructure network, some components of the state information in the original MDP model have implications on other components (if, for example, a node's fault is removed and an immediate successor node remains without power, then this eliminates uncertainty about whether or not the successor node contains a fault). Leveraging this, we introduce new state aggregation methods that significantly reduce the state space while maintaining competitive solution quality. We test our solution methods through extensive computational experiments based on information level, infrastructure network types and sizes, and repair duration, and compare our methods against benchmark solution methods.

The rest of this paper is organized as follows. In Section \ref{Sec:Lit}, we review related literature on the logistics of infrastructure network restoration. Sections \ref{Sec:Des} and \ref{Sec:Form} formally define the problem and provide a model formulation, respectively. In Section \ref{Sec:Approach}, we present our solution methodology. Section \ref{Sec:Comp} presents and discusses the results of a broad set of computational tests. Finally, Section \ref{Sec:Conc} summarizes our work and discusses potential further avenues for research.

\section{Literature Review}

\label{Sec:Lit}

The body of literature on the recovery of power networks after a disaster is relatively new. Strategies for power network recovery can be categorized into three groups: managing crew routing, portable energy source system deployment, and network reconfiguration by rerouting power (using the power flow equations). This work focuses on the former class of crew routing problems that occur in the last-mile of a distribution network (the network between substations and customers), whereas the latter two approaches focus on restoring flow throughout the transmission network (the network between power plants and substations). 
This distinction is crucial, as understanding the different network segments and their unique restoration is challenging. \cite{biswas2024review} provide an in-depth review of response strategies for restoring infrastructure networks in the aftermath of disasters. They systematically examine quantitative methods developed to facilitate the restoration of various infrastructure networks, categorizing these methods based on the type of network under consideration. The authors place special emphasis on resource allocation, scheduling, routing, and repair efforts, specifically focusing on networks related to power, roads, water, oil, and gas. \cite{salehipour2011efficient} and \cite{ccavdar2022repair} each study a single-crew routing problem for distribution network restoration under complete information on network fault locations. The former incorporates a meta-heuristic strategy that involves Greedy Randomized Adaptive Search Procedure (GRASP) and variable neighborhood search for determining a route, whereas the latter introduces an exact solution method that leverages a bi-directional dynamic programming approach along with structural results to partially characterize sub-optimal solutions. \cite{van2011vehicle}, \cite{morshedlou2018work} and \cite{morshedlou2021heuristic} study multiple-crew routing problems. The first of these aims to minimize the total outage time by generating precedence relationships and integrating them into a vehicle-routing sub-problem to enhance solution quality for large instances. On the other hand, \cite{morshedlou2021heuristic} propose a heuristic algorithm incorporating density-based mapping and solution clustering for crew allocation, whereas \cite{morshedlou2018work} focus on maximizing power network resilience over time by developing two mixed-integer programming (MIP) models, introducing binary and proportional active restoration crew routing models. The binary active restoration model operates on the principle that each component of the system can be in one of two states: active or inactive, whereas proportional active restoration considers the partial allocation of resources to an arc whose capacity is reduced due to disruptions. \cite{li2021hybrid} address an integrated repair crew routing problem in distribution network restoration accounting for multiple uncertainties, such as crew travel time, repair time, electricity demand, and photovoltaic (PV) generation. The authors employ a multi-stage model that includes routing, distribution network reconfiguration, and robust optimization stages. 
They utilize stochastic programming, robust optimization, column generation, and progressive hedging to minimize the total system operating and demand curtailment costs. 

In network reconfiguration using portable energy systems, \cite{lei2019resilient} simultaneously tackle repair crew and portable energy source decisions to restore service efficiently. They model the scheduling and routing of repair crews and energy systems within a distribution system to reconfigure the networks to form microgrids. 
Meanwhile, \cite{arif2018optimizing} extend this coordination to uncertain scenarios, introducing a two-stage stochastic mixed integer linear programming model 
with a progressive hedging algorithm, dispatching repair crews and reconfiguring 
the power network sequentially. \cite{arif2017power} employ a two-step strategy, initially clustering damaged components based on required travel distances and repair crew availability. Subsequently, they design a mixed-integer linear programming (MILP) model, simultaneously optimizing portable energy source system dispatching and repair crew routing decisions. \cite{ding2020multiperiod} formulate an MILP model to address the synchronization of mobile energy source systems, repair crews, and network microgrids for power restoration. They introduce an algorithm that maintains the solution quality of CPLEX but significantly enhances the computational speed with respect to CPLEX by adjusting the weights of the critical loads in the objective function. This is achieved by adding an auxiliary term that modifies the coefficients of specific integer variables, particularly those related to the status of power lines, which are critical in determining the system's configuration during restoration. 
\cite{arab2016electric} focus on optimal repair schedules and network reconfiguration plans considering the criticality of loads. They employ an MIP model and use Benders' decomposition to solve it. \cite{chen2018toward} develop an MILP model addressing crew routing, network reconfiguration, and the sequence in which electrical loads are restored such that the restoration process maximizes system stability, safety, and efficiency. 
\cite{carvalho2007dynamic} propose network reconfiguration operations with an optimal sequence determined using dynamic programming, employing methods to reduce state complexity. \cite{van2015transmission} present a two-stage approach, utilizing an MIP model for power flow equations in the first stage and extensive neighborhood search for repair crew routing in the second. \cite{nurre2012restoring} integrate network design and crew scheduling, formulating an integer programming model with valid inequalities and employing a rule-based dispatching heuristic for efficient restoration.

Our underlying problem is a crew routing problem under incomplete information on the status of network nodes, with the objective of minimizing service disruption time.  
A small number of prior studies consider related problems where information is gradually revealed over a rolling planning horizon, and while these methodologies are closely related to ours, they do not address the critical precedence relationships that exist in infrastructure networks. 
For example, \cite{ulmer2016rollout,ulmer2018a,ulmer2018b,ulmer2019offline}, and \cite{ulmer2019anticipation} focus on optimizing vehicle routing to maximize the number of people served or deliveries completed within limited hours. These studies analyze vehicle routes where information about customer requests is gradually revealed during the planning horizon, but they do not consider the precedence constraints that are crucial in infrastructure networks. 
This stream of work relates to infrastructure network restoration by efficiently allocating limited resources (e.g., service vehicles or road maintenance crews) to maximize restoration within a constrained time frame, as real-time information about the status of nodes in the power network is revealed. 
Similarly, \cite{ulmer2020modeling} propose a route-based Markov decision process (MDP) that extends conventional MDP-based vehicle routing approaches by incorporating problem-specific route plans into the set of feasible actions. Problem-specific route plans correspond to predefined sequences of actions that are tailored to the unique characteristics of a given problem. These plans take into account factors such as the timing of pickups and deliveries, the availability of service vehicles, and the need to adhere to precedence constraints (e.g., ensuring that a pickup occurs before a delivery). 
They employ dynamic programming principles to solve the corresponding MDP. \cite{ulusan2021approximate} study a network recovery problem under uncertainty but do not consider the real-time status of the nodes or the inter-dependencies between different networks. 
They frame the problem as an MDP and apply approximate dynamic programming techniques to derive solutions that are close to optimal. Utilizing an approximation approach, they estimate the value functions by representing them as linear combinations of states. \cite{goldbeck2020optimal} propose a multi-stage stochastic programming model that jointly optimizes pre-disruption investments in network capacities and repair capabilities, as well as post-disruption operational adjustments and the allocation of repair resources. The model is dynamic, considering multiple time steps and stages of decision-making, which permits capturing the progression of disruptions and the effectiveness of recovery strategies over time. 

The existing literature on repair crew routing for infrastructure networks is still very limited, particularly in terms of understanding of the impact of precedence relations and incomplete information on the ability to efficiently determine high-quality solutions. 
Our study contributes to this under-explored area by introducing a new solution approach for the single crew routing problem under incomplete information with an objective of minimizing total service disruption time. The following section introduces the problem, provides insights into the specific challenges, constraints, and objectives, and presents a problem formulation. 

\section{Problem Definition}
\label{Sec:Des}

We consider a disrupted power distribution network where the planner knows which locations have lost power but not necessarily which locations actually have a fault and require repair. 
Our goal is to determine a route for the repair crew that determines the order in which network locations should be visited for inspection and repair in order to minimize total service disruption time. 
Before going into the problem description, we first provide some background on the context.

\begin{figure}[h]
  \centering \includegraphics[width=0.35\linewidth]{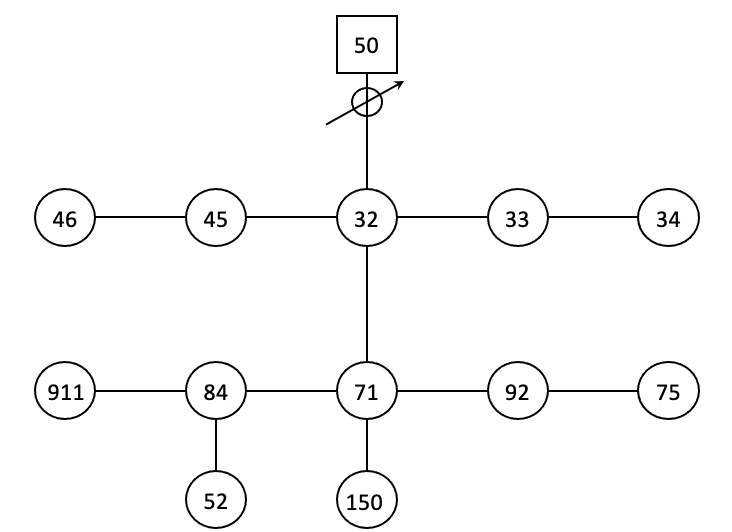}
  \caption{A representation of a distribution network (IEEE 13-node instance).}
  \label{Fig:13feeder}
\end{figure}

We consider two networks, namely the disrupted infrastructure network and the road network that permits travel to nodes on the infrastructure network. Although we focus on a power network as the infrastructure network, our results are generalizable for other infrastructure networks. 
Power distribution networks commonly have a directed-tree structure, i.e., there are no cycles.
Figure \ref{Fig:13feeder} displays the IEEE 13-feeder instance as an example distribution network. 
In the figure, node 50 is the source, and others are demand nodes that receive flow from the source.
If, for example, there is a fault at node 84, then in addition to node 84, its successors (i.e., nodes 52 and 911) also lose service due to the precedence relationship determined by the direction of power flow from the source. 
In the case that neither node 52, 84, nor 911 has service, then we can immediately infer that node 84 has a fault. 
However, this information does not reveal whether nodes 52 and 911 have faults as well because the precedence relationships on the network hide this information.
Hence, if the repair crew visits node 52, upon physical inspection, they may observe that there is no fault, resulting in an unnecessary trip. In case of incomplete information on the actual status of node faults, unnecessary trips may not be entirely avoidable. 

We denote the underlying infrastructure network by $G_{\mathcal{I}}=(N_{\mathcal{I}}, A_{\mathcal{I}})$, where $N_{\mathcal{I}}$ is the set of nodes on the distribution network and $A_{\mathcal{I}}$ is the set of arcs connecting these nodes. $N_{\mathcal{I}}=\{1, 2, \dots, n\}$ where node 1 represents the single source, and other nodes represent demand locations.
The distribution network is a directed tree rooted at node 1. The set of arcs $A_{\mathcal{I}}=\{(i,j): i, j \in N_{\mathcal{I}} \text{ and } i\neq j\}$ represents the immediate precedence relationships on the power network. An arc $(i,j)$ implies that node $j$ is an immediate successor of node $i$ on the distribution network. Since the distribution network is a tree, there is a unique path between the source node and each demand node $i$ denoted by $P_{1i}$ that also includes node $i$ itself.
For a path $P_{1i}$, any node $j$ on the path is a predecessor of node $i$, which we denote by $i \preceq j$, implying that node $i$ cannot receive service before node $j$. 
The service at node $i$ can only be recovered if all the faults on the path $P_{1i}$ are removed. 
If a node has a fault, we assume that it takes a duration of $s$ time units for the crew to repair.

The road network is denoted by $G=(N, A)$ where $N$ is the set of nodes on the road network and $A$ is the set of edges connecting these nodes.
We assume that the road network is an undirected complete graph.
Without loss of generality, we can assume that $N=N_{\mathcal{I}} \cup \{0\}$ where node 0 is the depot where the repair crew is located. 
We define $d_{ij}$ as the travel time on the shortest path between nodes $i$ and $j$ in $N$, which is measured in the same units as the repair time.

We consider a disrupted distribution network $G_{\mathcal{I}}=(N_{\mathcal{I}},A_{\mathcal{I}})$ where a subset of nodes experiences a power service loss, and the repair crew does not necessarily know the locations of all faults, except for those where it is immediately observable prior to physical inspection (e.g., where a node's immediate predecessor has power). 
Assuming that all nodes are subject to the same level of disruption, we denote the probability of a node having a fault by $p$ if the true status cannot be directly observed.
We present a summary of our notation in Table \ref{Tab:Notation}.

Our goal is to determine a crew route to visit and repair all faults on the infrastructure network with minimum expected total service disruption time. Considering the dynamic nature of the decisions due to the change in available information over time, we formulate this problem as a discrete-time, finite-horizon Markov decision process (MDP). We present the elements of our model and the formulation in the next section.

\section{Formulation}

\label{Sec:Form}
In our MDP formulation, the decision epochs represent the points in time when the crew must determine its next location. 
Thus, such a decision epoch occurs after repairing an existing fault or observing that the current node has no fault, and then updating node status information accordingly. 
Hence, each decision epoch corresponds to a departure from a node on the infrastructure network. 
The state is defined based on information on the crew's location and the status of all nodes in the infrastructure network. 
We denote the state at decision epoch $t$ by $S_t=(L, U^+, U_0^-, U_1^-, U_p^-)$ where $L \in N$ is the current location of the crew, $U^+$ is the set of nodes that have service, and $U_q^-$ is the set of nodes without service and faulty with probability $q$, with $q\in \{0, p, 1\}$ and $0<p<1$. The subsets $U^+$, $U_0^-$, $U_1^-$ and $U_p^-$ are mutually exclusive and collectively exhaustive for $N_{\mathcal{I}}$.
A feasible state must have the following properties:
\begin{itemize}
\item If node $i \in U^+$, every node $j\in P_{1i}$ must also be in $U^+$.
\item If node $i \in U_0^-$, there is at least one node $j \in P_{1i}$ that is not in $U^+$. 
\item If node $i \in U_1^-$, every node $j \in P_{1i}$ and $j\neq i$ must be in $U^+$.
\item If $i \in U_p^-$, there is at least one node $j \in P_{1i}$ and $j\neq i$ that is not in $U^+$.
\item If node $i$ is in $U_0^-$, then this implies that node $i$ was already visited by the repair crew.  However, a node $i$ that is in $U_1^- \cup U_p^-$ has not yet been visited. Hence, $L \in U^+ \cup U_0^-$.
\end{itemize}

We denote the set of all feasible states by $\mathcal{S}$.  Any state with $U^+=N_{\mathcal{I}}$ implies that service restoration has been completed.

\begin{table}[h]
\begin{center}
\caption{Summary of notation.}
 \scalefont{0.8}
\begin{tabular}{c|p{11cm}}

    \hline Notation & Description \\\hline
    $i$, $j$ & node indices \\
    $t$ & index for decision epochs\\
    $G_{\mathcal{I}}=(N_{\mathcal{I}}, A_{\mathcal{I}})$  & infrastructure network \\
    $G=(N, A)$ & road network\\
    $P_{1i}$ & the unique path between the source and node $i$ on $G_{\mathcal{I}}=(V_{\mathcal{I}}, A_{\mathcal{I}})$ \\ 
    $d_{ij}$  & travel time over edge $(i,j)$\\ 
    $s$ & service time at fault locations\\
    $p$ & node fault probability under incomplete information\\
    \hline
    $U^{+}$ & set of nodes with service \\
    $U^{-}_1$ & set of nodes without service and faulty with probability 1 \\  
    $U^{-}_0$ & set of nodes without service and faulty with probability 0 \\  
    $U^{-}_p$ & set of nodes without service and faulty with probability $p$\\
    $L$ & current location of the crew\\ 
    $S_t=(L, U^+, U^{-}_0,U^{-}_1,U^{-}_p)$ & state at decision epoch $t$\\        
    $\mathcal{S}$ & set of all possible states \\
    $a_t$ & action at decision epoch $t$, i.e., the next node to go to \\
    $r_j$ & recovery time of node $j$\\
    $c(S_t,a_t)$ & cost of action $a_t$ in state $S$: additional expected service disruption time for all nodes that occurs between the time the crew leaves the current node and finalizes the repair at node $a_t$ (if fault exists). \\
    \hline
\end{tabular}
\label{Tab:Notation}
\end{center}
\end{table}

\textbf{Actions and one-stage cost: }
In each state where there is at least one remaining node without service, we need to determine the crew's next location. We denote the corresponding action by $a_t$, where $a_t \in U_1^- \cup U_p^-$. Given a state $S_t=(L, U^+, U^{-}_0,U^{-}_1,U^{-}_p)$, the cost of action $a_t$ is the expected increase in the total service disruption time for all nodes in the infrastructure network between the time the crew departs node $L$ and completes the repair at node $a$ (if there is a fault at node $a_t$). This cost is the product of the travel plus repair time and the total number of nodes without service during that time.
We denote this cost by $c(S_t,a_t)$, which is calculated as 
\begin{eqnarray}\label{Eq:Cost}
    c(S_t,a_t) = \begin{cases}
      \big(d_{La_t} + s \big)\sum_{j \in  V_{\mathcal{I}}\setminus U^{+}}   \mathbbm{1}_{j }  , & \text{if $a_t \in U^{-}_1$}, \\
      (d_{La_t} +p s)\sum_{j \in V_{\mathcal{I}}\setminus U^{+}} \mathbbm{1}_{j }, & \text{if $a_t \in U^{-}_p$}. \\
\end{cases} 
\end{eqnarray}
where $\mathbbm{1}_{j }$
is a binary variable equal to 1 if node $j$ belongs to the set of nodes without power, and 0 otherwise. The sum of these binary indicator variables is the number of nodes without power. Notice that if $a_t\in U_1^-$, then there is a fault at node $a_t$ that takes $s$ time units to repair. On the other hand, if $a_t\in U_p^-$, then the crew will need to perform a repair with probability $p$.

\textbf{Transition probabilities: } With each crew action, new information is revealed that changes the status of a subset of nodes. When the crew visits a node $a\in U_p^-$, then either a fault is observed and repaired, or no fault exists. In both cases, node $a_t$ is labeled as without fault, but cannot gain service until all faults of predecessor nodes are repaired.
On the other hand, if the crew visits a node $a_t \in U_1^-$, the fault at this location is repaired, and node $a_t$ gains service. In addition, this reveals new information about the successors of node $a_t$. This information depends on the true status of such nodes (i.e., faulty or not) and the previous crew actions. 

\begin{figure}[h]
\begin{center}
\begin{tikzpicture}[scale=0.5][>=stealth]
\node[el1] (e1)at (12,5) {\color{black}$U^{+}$};
\node[el1] (e2)at (9,5) {\color{black}$U^{-}_0$};
\node[el1] (e3)at (6,5) {\color{black}$U^{-}_1$};
\node[el1] (e4)at (3,5) {\color{black}$U^{-}_p$};
\node[el1] (e5)at (12,0) {\color{black}$\bar{U}^{+}$};
\node[el1] (e6)at (9,0) {\color{black}$\bar{U}^{-}_0$};
\node[el1] (e7)at (6,0) {\color{black}$\bar{U}^{-}_1$};
\node[el1] (e8)at (3,0) {\color{black}$\bar{U}^{-}_p$};
\draw [black,line width=0.5mm][ ->] (e1) to node[right]{\color{red}} (e5);
\draw [black,line width=0.9mm][dashed, ->] (e4) to [left] node [above right]{\color{red}} (e6);
\draw [black,line width=0.4mm][dashed, ->] (e2) to node[right]{\color{red}} (e6);
\draw [black,line width=0.9mm][->] (e3) to node[right]{\color{red}} (e7);
\draw [black,line width=0.9mm][dashed, ->] (e4) to node[right]{\color{red}} (e8);
\end{tikzpicture}
\qquad \qquad 
\begin{tikzpicture}[scale=0.5][>=stealth]
\node[el1] (e1)at (12,5) {\color{black}$U^{+}$};
\node[el1] (e2)at (9,5) {\color{black}$U^{-}_0$};
\node[el1] (e3)at (6,5) {\color{black}$U^{-}_1$};
\node[el1] (e4)at (3,5) {\color{black}$U^{-}_p$};
\node[el1] (e5)at (12,0) {\color{black}$\bar{U}^{+}$};
\node[el1] (e6)at (9,0) {\color{black}$\bar{U}^{-}_0$};
\node[el1] (e7)at (6,0) {\color{black}$\bar{U}^{-}_1$};
\node[el1] (e8)at (3,0) {\color{black}$\bar{U}^{-}_p$};
\draw [black,line width=0.5mm][->] (e1) to [left] node [above right]{\color{red}} (e5);
\draw [black,line width=0.9mm][dashed, ->] (e4) to node[right]{\color{red}} (e7);
\draw [black,line width=0.9mm][dashed, ->] (e4) to node[right]{\color{red}} (e5);
\draw [black,line width=0.9mm][ ->] (e3) to [left] node [above right]{\color{red}} (e5);
\draw [black,line width=0.4mm][dashed, ->]  (e2) to [left] node [above right]{\color{red}} (e5);
\draw [black,line width=0.4mm][dashed, ->]  (e2) to node[right]{\color{red}} (e6);
\draw [black,line width=0.9mm][ ->] (e3) to node[right]{\color{red}} (e7);
\draw [black,line width=0.9mm][dashed, ->] (e4) to node[right]{\color{red}} (e8);
\end{tikzpicture}
\end{center}
\caption{Possible node transitions  for $a \in U^{-}_p$ (left) and for $a \in U^{-}_1$ (right).} \label{fig:NodeTransition}
\end{figure}
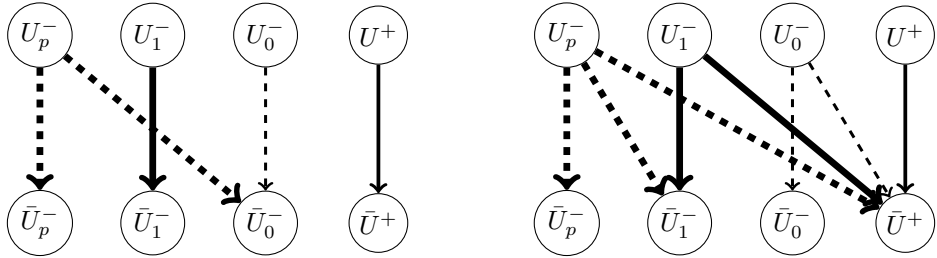

The transition probability is denoted by $P(S_{t+1}|S_t, a_t)$, which corresponds to the probability of the state transition  
$S_t = (L,U^{+}, U^{-}_1, U^{-}_0, U^{-}_p) \rightarrow S_{t+1} = (\bar{L}, \bar{U}^{+},\bar{U}^{-}_1,\bar{U}^{-}_0, \bar{U}^{-}_p)$ when action $a_t$ is taken in state $S_t$. In Figure \ref{fig:NodeTransition}, we represent all possible node transitions that may result from different actions, depending on whether the crew travels to a faulty node $a_t\in U^-_1$ or $a_t\in U^-_p$.

Computing the transition probability when $a_t\in U^-_p$ is simpler than when $a_t\in U^-_1$. When $a_t\in U^-_p$, the state transition is deterministic; node $a_t$ moves from $U^-_p$ to $\bar{U}^-_0$,  the current location is updated to $a_t$, and there is no change in the status of other nodes. On the other hand, when $a_t\in U^-_1$, we know that node $a_t$ is the first node on path $P_{1a}$ with a fault, and repairing it will restore service at this location. Therefore, node $a_t$ will move to $\bar{U}^+$ and will be the current node. In addition, some of the successors of node $a_t$ may also gain service, depending on the true status of the nodes and network structure. Hence, this part of the transition is probabilistic. We therefore calculate the transition probabilities as 
\begin{equation}\label{Eq:TransitionProb}
P(S_{t+1}|S_t, a) = 
\begin{cases}
  \Big((1 - p)^{\lvert \bar{U}^{+} \rvert + \lvert \bar{U}^{-}_0 \rvert - \lvert U^{+} \rvert - \lvert U^{-}_0 \rvert - 1} \cdot p^{\lvert \bar{U}^{-}_1 \rvert - \lvert U^{-}_1 \rvert + 1}\Big) \mathbbm{1}_{S_t, S_{t+1}, a}, & \text{if } a \in U^{-}_1 \\
  \mathbbm{1}_{S_t, S_{t+1}, a},  & \text{if } a \in U^{-}_p,
\end{cases}
\end{equation}
where $\mathbbm{1}_{S_t, S_{t+1}, a}$ is a binary indicator for whether the deterministic node transitions between states $S_t$ and $S_{t+1}$ as a result of action $a_t$ has occurred (1) or not (0). The quantity $\lvert \bar{U}^{+} \rvert + \lvert \bar{U}^{-}_0 \rvert - \lvert U^{+} \rvert - \lvert U^{-}_0 \rvert - 1$ determines the number of nodes observed to be not faulty, while $\lvert \bar{U}^{-}_1 \rvert - \lvert U^{-}_1 \rvert + 1$ is the number of nodes observed to be faulty after taking action $a_t\in U^-_1$. We present pseudo-code to calculate the transition probabilities in Algorithm \ref{Alg:TransitionPseudo} in Appendix \ref{Ap:TransProp}.

\textbf{Objective function:  }Our goal is to minimize the total expected service disruption time. Given a state $S_t$, we denote the minimum total expected service disruption time until completion by $H_t(S_t)$, where $t$ is the index for decision epochs, $t=0, 1, \dots, n$.
Combining the cost function in Equation (\ref{Eq:Cost}) and the transition probabilities in Equation (\ref{Eq:TransitionProb}), we can write the Bellman equation for minimum expected service disruption time as:
\begin{equation} \label{Eq:Obj}
\begin{split}
H_t (S_t) &= \min_{a\in \mathcal{A}_{S_{t}}} \Big\{c(S_t, a) + \sum_{S_{t+1} \in \mathcal{S}} P(S_{t+1}|S_t, a) H_{t+1}(S_{t+1})\Big\}  \\
 & = \min_{a \in \mathcal{A}_{S_t}} \Big\{c(S_t, a) + \displaystyle \mathop{\mathbb{E}}\{ H_{t+1}(S_{t+1})|S_t, a\}\Big\} ,
\end{split}
\end{equation}
where $\mathcal{A}_{S_t}$ is the set of available actions in state $S_t$, with the termination condition
\begin{equation}\label{bellman}
H_t ( S_t=(L, N_I, \emptyset, \emptyset, \emptyset)) =0, \quad \forall L \in N_I.
\end{equation}
 We call this problem the Traveling Repairman Network Restoration Problem (TRNRP). Solving the TRNRP presents significant computational challenges. Some of these challenges are inherent in classical routing decisions, where the problem size increases exponentially in the number of nodes. For the TRNRP, incomplete information and the dependence on the infrastructure network further complicate the state transitions. 
To solve the problem efficiently, we use an approximate dynamic programming approach that is strengthened by structural results. We also propose state aggregation methods to reduce the search space. The following section presents our solution approach. 

\section{Solution Approach}
\label{Sec:Approach}
Our solution approach determines the repair routes with the minimum total expected service disruption time using dynamic programming (DP). However, due to the exponential nature of the TRNRP state space, the exact implementation of the DP imposes computational limitations. 
To address these challenges, we develop an approximate DP (ADP) approach where we perform value approximation through Reinforcement Learning (RL). Our RL application follows \cite{powell2007approximate}.
We also introduce some structural results to identify suboptimal routes. In addition, we present problem-specific state aggregation methods to increase computational efficiency. We present the details of our solution approach below.

\textbf{Value approximation using RL: }
Recall that the state $S$ in our MDP formulation represents the information available to the repair crew just before taking a new action. Thus, our original states are called \emph{pre-decision states} \citep{powell2007approximate}. 
To enable value approximation via RL, we introduce \emph{post-decision states}. Post-decision states store the same type of information as in the pre-decision states, incorporating the most recent action taken, but without new information on which nodes will have gained power as a result of the chosen action.
We denote the post-decision state that occurs after taking action $a_t$ on a pre-decision state $S_t$ by $S^a_t$, and let $H^a_t(S^a_t)$ denote the value of the post-decision state.  
Let $g(\cdot, \cdot)$ be a state transformation function that takes the current pre-decision state ($S_t$) or the post-decision state ($S^a_t)$ as its first argument and information $\omega$ revealing updates on the status of all nodes as the second argument. The function $g(\cdot, \cdot)$ determines the new state, i.e., the post-decision state that follows $S_t$ or the next pre-decision state after $S^a_t$. 
The information contained in $\omega$ can be one of the following:
\begin{itemize}
\item[(i)] the current node being repaired as a result of action $a_t$, but no new information on the subsequent condition of the successors of node $a_t$. In this case, $g(S_t, \omega)$ provides the post-decision state $S^a_t$.
\item[(ii)] information on which of the successors of node $a_t$ have gained service as a result of the action. In this case $g(S^a_t, \omega)$ provides the new pre-decision state $S_{t+1}$. 
\end{itemize}

The relation between the values of $H_t(S_t)$ and $H^a_t(S^a_t)$ can be written as follows:
\begin{equation}\label{Eq:PostValue}
    H_{t}(S_t) = \min_{a_t \in \mathcal{A}_{S_t}} \{c(S_t, a_t) + H^a_t(S^a_t)\}
\end{equation}
\begin{equation}\label{Eq:PostExpected}
    H^a_{t}(S^a_t) = \displaystyle \mathop{\mathbb{E}}[H_{t+1}(S_{t+1})| S^a_t]
\end{equation}
where  $S^a_t = g(S_t, \omega_1)$ and $\omega_1$ includes the information that node $a_t$ is being repaired (or visited), and $S_{t+1} = g(S^{a}_t, \omega_2)$ and $\omega_2$ includes information of which set of new nodes gain service after repairing node $a_t$ in state $S_t$.
 
By substituting Equation \eqref{Eq:PostExpected} in Equation \eqref{Eq:PostValue}, we obtain Equation \eqref{Eq:Obj}, which is the original problem we wish to solve. Alternatively, by substituting Equation \eqref{Eq:PostValue} in Equation \eqref{Eq:PostExpected}, we obtain the Bellman Equation \eqref{Eq:ExpectedMin} in terms of the post-decision state variables. 
Solving Equation \eqref{Eq:ExpectedMin} by computing the expected value of the minimization is faster compared to solving Equation \eqref{Eq:Obj}, where we minimize the expected value of being in state $S_t$. Equation \eqref{Eq:Obj} requires minimizing the expected value by considering all possible state transitions for each possible action, leading to a high computational burden. 
In contrast, in Equation \eqref{Eq:ExpectedMin}, the minimization problem inside the expected value is a deterministic optimization problem if the value of $H^a_{t+1}(S^a_{t+1})$ is known for all $a_{t+1} \in \mathcal{A}_{S_{t+1}}$.
\begin{equation}\label{Eq:ExpectedMin}
    H^a_{t}(S^a_t) = \displaystyle \mathop{\mathbb{E}} \bigg{[}\underbrace{\min_{a_{t+1} \in \mathcal{A}_{S_{t+1}}} \{c(S_{t+1}, a_{t+1}) + H^a_{t+1}(S^a_{t+1})\} }_{H_{t+1}(S_{t+1})}| S^a_t \bigg{]}
\end{equation}

Although we do not know $H^a_{t+1}(S^a_{t+1})$, we can approximate its value using simulation, where we let $\hat{H_t}^a(S^a_t)$ denote this approximate value. During the simulation, we solve the TRNRP repeatedly based on randomly generated instances. We refer to the solution of any given instance until completion an \emph{iteration}. 
In each iteration, decisions are made based on current approximate values, and a new approximation is computed for the state under consideration based on the results of the iteration. 
Using this new value, the approximate value is updated. 
This process allows value propagation between reachable states. 
That is, when we are choosing an action in a state, we consider the immediate cost of the action and the expected value of reachable states. Therefore, as the simulation continues, any updates to a state's value affect the values of all states from which the state is reachable, which leads to value propagation.
The simulation ends when a set of stopping criteria (which we will describe shortly) are met. We describe the details of the approximation below.

Let $\hat{H}_t^{a,n}(S^{a}_t)$ be the approximate value of post-decision state $S^a_t$ after $n$ iterations of the simulation. The value $\hat{H}_t^{a,n}(S^{a}_t)$ incorporates the simulation results for iterations $1, 2, \dots, n$. In iteration $n$, when determining which action to take in state $S_t$, we consider the value of $\hat{H}_t^{a,n-1}(S^{a}_t)$ for all possible actions in the current state that lead to a corresponding post-decision state $S^{a}_t$, and choose the action with minimum total expected service disruption time.
Once an action is selected, we calculate a new value for the specific post-decision state associated with the chosen action. 
Let $\hat{h}_t^n (S_t^{a, n})$ be the sample realization of the post-decision state $S_t^{a, n}$ in the $n^{th}$ iteration, i.e., 
\begin{align}\label{min_postdecision}
\hat{h}^n_{t}(S_t^{a, n})  =& \min_{a_t \in \mathcal{A}_{S_t}} \{c(S_t, a_t) + \hat{H}^{a,n-1}_{t}(g(S_t, \omega(a_t))\}\\
=& \min_{a_t \in \mathcal{A}_{S_t}} \{c(S_t, a_t) + \hat{H}^{a,n-1}_{t}(S^{a}_t)\},
\end{align}

where $g(S_t, \omega(a_t))$ determines the post-decision state based on the new information $\omega(a_t)$ that arises from the action under consideration.
For each new $\hat{h}^n_t (S_t^{a, n})$, we use exponential smoothing to compute an updated value of the corresponding post-decision state as 
\begin{equation}
\hat{H}^{a,n}_{t}(S^{a,n}_t) = (1-\alpha_{n-1})\hat{H}^{a,n-1}_{t}(S^{a,n}_t) + \alpha_{n-1} \hat{h}^n_{t}(S_t^{a, n}),
\end{equation}
where $\alpha_{n-1}$ is the smoothing factor at iteration $n$. During the simulation, as the frequency of observing a state increases, we reduce the smoothing constant. Accordingly, $\alpha_{n-1}$ is computed as  
\begin{equation*}
    \alpha_{n-1} = \frac{1}{N(S^{a,n}_t)},
\end{equation*}
where $N(S^{a,n}_t)$ is the number of times the post-decision state $S^a_t$ has been observed over the $n$ iterations. In addition to this update mechanism, to avoid stalling at a local minimum during simulation, we use a variant in which we choose an action randomly, rather than minimizing the expected cost-to-go. The choice of a random decision is known as \emph{exploration}, whereas decisions based on the minimum expected cost-to-go are referred to as \emph{exploitation}. 
To determine when to use exploration versus exploitation, we use an $\epsilon$-greedy method, where $\epsilon$ corresponds to the fraction of decisions where we choose actions based on exploration.  We explore at a higher rate in earlier iterations, reducing exploration as the simulation progresses. Accordingly, the value of $\epsilon$ in iteration $n$ is calculated as
\begin{equation*}
    \epsilon^{n} = \frac{\alpha}{N(S^{a,n}_t)}.
\end{equation*}

We start the simulation with a randomly generated fault realization for the network nodes, wherein all nodes have a service loss and a subset of nodes contains faults, although the dispatcher does not know the actual fault locations. 
During the simulation, after each action, we update the state information based on the true nature of the instance.

The simulation continues until the following stopping conditions are met.
After 100,000 simulation iterations, we begin checking the progression of approximate state values every 10,000 iterations, which we call a batch of iterations.
If, in the most recent batch, we have observed a new state, then we continue with another batch of 10,000 iterations.  On the other hand, if the set of observed states remains the same between batches, we consider the value change in each of the most frequently observed states, defined as those states that have been visited more in more than 75\% of all prior iterations. That is, for each of the most frequent states, we compute the absolute change in the approximate value, i.e, $|\hat{H}^{a,n+10,000}_t (S^{a,n+10,000}_t) - \hat{H}^{a,n}_t (S^{a,n}_t)|$, in the last batch.
Then, we compute the maximum absolute change over the set of the most frequently observed states. Let us denote this value for batch $k$ by $\Delta_k$.
If the maximum of the last three recorded absolute differences, i.e., $\max \{\Delta_{k-2}, \Delta_{k-1}, \Delta_k\} $, is less than a predefined threshold $\delta$, we terminate the simulation. The simulation then returns a look-up table including the approximate values of the post-decision states. Mathematically, this condition is represented as
 \scalefont{0.85}
\begin{equation}\label{stopping}
\max_{j \in \{k-2, k-1, k\}} 
\Bigg( 
    \max_{S \in \mathcal{S}_{\text{frequent}}} 
    \Big| 
        \hat{H}^{a, n + 10,000(k+1-j)}_t\left(S^{a, n + 10,000(k+1-j)}_t\right) \nonumber 
        - 
        \hat{H}^{a, n + 10,000(k-j)}_t\left(S^{a, n + 10,000(k-j)}_t\right) 
    \Big| 
\Bigg) < \gamma.
\end{equation}

\normalsize

We present a pseudo-code for our algorithm for computing the approximate look-up table in Algorithm \ref{alg:cap1} in Appendix \ref{Alg:TransitionPseudo}. During route planning, these values are used to choose the action with minimum cost-to-go in real time. We call this solution method the Network Restoration Routing ($NRR$) approximation algorithm.

\subsection{Refinements}
The preceding section introduced the $NRR$ approximation algorithm designed for determining a repair crew route under incomplete information on node status. During the simulation, particularly in early stages where the quality of approximations is low, sub-optimal actions may be selected, which increases the required number of iterations until completion. To mitigate this effect, we next propose some refinements that rely on identifying a subset of sub-optimal actions to avoid.

Let $\mathcal{P}$ denote a feasible offline solution for the crew routing problem and let $i$ be a node on the associated path corresponding to the crew route. Since $\mathcal{P}$ is an offline solution, some of the nodes on the route may gain service prior to visiting the node, in which case the node can be skipped on the path. 

\begin{proposition}
    Consider nodes ${i, j , k} \in N_{\mathcal{I}}$ such that $k \preceq j$ on the infrastructure network. If $d_{ij} \leq d_{ik}$ and $d_{kh} \leq d_{jh}$ and node $h$ is not the depot, or  $d_{ij} \leq d_{ik}$  and node $h$ is the depot, then any offline solution to the crew routing problem that has $i\rightarrow j \rightarrow k \rightarrow h$ is sub-optimal and therefore can be eliminated from consideration.
    \label{Prop:PrereqRefinement}
\end{proposition}

The proof is in Appendix \ref{Proof_of_Proposition_1}. Proposition 1 implies that a node should not be visited before its predecessor on the infrastructure network unless doing so would reduce the travel time. Otherwise, the crew incurs greater travel distance and time with less information on the network's fault status and a greater number of nodes without service. Motivated by this result, we next define some additional terms needed to obtain more general results.

For an offline path $\mathcal{P}$, we denote $Q_{q}^{\mathcal{P}}$ as the total number of nodes without service at the time the crew leaves the $q^{th}$ node on the path, and $d_q^{\mathcal{P}}$ as the travel time on the $q^{th}$ arc on path $\mathcal{P}$.
\begin{definition}
     If $d_q^{\mathcal{P}}$ cannot be decreased using a $k$-opt local search for any $q$, we call this path $k$-opt guaranteed with respect to travel time.
\end{definition}

\begin{definition}
If $Q_{q}^{\mathcal{P}}$ cannot be decreased for any $q$ using a $k$-opt local search, we call this path $k$-opt guaranteed with respect to information.
\end{definition}
\begin{definition}
If a path  $\mathcal{P}$ is both $k$-opt guaranteed with respect to travel time and information, we refer to such a path as double $k$-opt guaranteed.
\end{definition}
\begin{proposition}
An optimal solution to the crew routing problem is double $k$-opt guaranteed for all  $k \geq 2$.
\label{Prop:DoubleSwap}
\end{proposition}

The proof is in Appendix \ref{Proof_of_Proposition_2}.

We call the solution method that uses these structural results in the decision process during the simulation the \emph{Strengthened Network Restoration Routing} ($sNRR$) approximation method. Since the $sNRR$ algorithm removes a set of sub-optimal actions at each epoch, it can facilitate faster computation of the approximate state values (i.e., training of the RL model).
Note that these structural results apply to offline solutions. However, we generate a route in an online fashion. Therefore, during the route generation process, we have limited information regarding the final route. Hence, when we incorporate these structural results into our original $NRR$ approximation method in a limited way, we lose some information.
Despite this, we observe that these refinements provide significant improvements during the training phase. 
In our preliminary numerical study, for experimental validation, we conduct simulations on 20 instances, obtained by varying the fault probability $p$ and the restoration duration $s$ for 16 and 18-node instances. 
To perform the comparison, in the first step, we run the original $NRR$ algorithm on these instances, letting each simulation run 4,000,000 iterations. Then, we determine the states that are visited at least 2.5\% of the time so that we account for commonly visited states. We randomly select five states out of these. In the second step, we run the $sNRR$ algorithm on the same instances, recording the iterations where we achieve convergence for the values of the corresponding states. Then, we compare the iterations where convergence is achieved for these states under $NRR$ and $sNRR$. We observe that over these 20 random instances, the refinements reduce the number of iterations to convergence by 52\% for the most commonly observed five states. Based on this comparison of $sNRR$ and $NRR$, we conclude that the former offers a significant computational advantage. Therefore, we use $sNRR$ as our base solution method.

Although these refinements allow us to execute the underlying approximate dynamic program faster with no compromise in solution quality, a challenge still remains due to the large size of the state space, which can preclude solving the problem for large instances. To address this challenge, in the next subsection, we introduce a set of state-space aggregation methods.

\subsection{State Aggregation}

Solving the TRNRP for large instances is a challenging task due to substantial memory requirements and the exponential growth of the state space. 
In our original formulation, we model the states to account for the exact information available at each decision epoch. However, the precedence relationships in the infrastructure network permit inferring partial information on a state based on other information elements of the state. Hence, the state space in our original model carries some unnecessary information. Therefore, to increase computational efficiency, we use state aggregation methods to reduce the size of the state space by grouping individual nodes with common types, but at the cost of losing precise information about the status of each individual node in the state space.

Considering the trade-off between solution quality and computation time, we present three aggregation methods for representing post-decision states. The aggregated states under these three methods are as follows.
\begin{itemize}
\item $SA^1: (L, |U^{+}|, U^{-}_1, |U^{-}_p|, |U^{-}_0|)$. Under this state aggregation approach, we keep track of the subset $ U^{-}_1$ as is and only store the cardinality of other subsets. This aggregation method is based on the fact that the repair crew must visit nodes in set $U^{-}_1$ during the repair process.
\item $SA^2: (L, |U^{+}|, |U^{-}_1|, |U^{-}_p|, |U^{-}_0|)$. This aggregation approach is less granular than the previous one, as it only stores the cardinality of all node subsets, resulting in a greater loss of specific information about individual nodes.
\item $SA^3: (L, |U^{+}|, |U^{-}_1|)$. Under the third aggregation approach, we only record the cardinality of subsets $ U^{+}$ and $U^{-}_1$, since node status does not change if the repair crew visits nodes in subsets $U^{-}_p$ and $U^{-}_0$.
\end{itemize}
State aggregation allows us to reduce the size of the lookup table, increasing the size of the instances we can solve. On the other hand, an aggregated state does not necessarily provide the full information needed to make an optimal decision on which node to visit next. Therefore, when we construct the lookup table using only including aggregated states, we keep track of the current pre-decision and post-decision states based on the realized information during the simulation as described earlier. Each post-decision state maps to an aggregated state with the corresponding approximate value, and this value is used to determine which action to take in the current state.
We present a flowchart depicting an iteration in our simulation that computes the approximate value of each state in creating the lookup table $\mathcal{L}$ in Figure \ref{fig:flowchart1} in Appendix \ref{App:FlowCharts}. Figure \ref{fig:flowchart2} in Appendix \ref{App:FlowCharts} shows the steps to determine the repair crew route for a given realization and lookup table $\mathcal{L}$.

The next section presents the results of a broad set of computational experiments to evaluate the performance of our proposed methods across different scenarios.

\section{Computational Experiments}
\label{Sec:Comp}
In this section, we present and discuss computational results from the application of our solution methods and compare them with benchmark heuristics. We use priority sequence ($PS$) and nearest neighbor ($NN$) as benchmark heuristics. In the $PS$ heuristic, we consider the number of successors of a node as its priority measure and, at each step, choose the node with the highest priority to visit next. The latter $NN$ heuristic follows a myopic policy based on distance and chooses the closest node that has not yet been visited at each step.

To test our solution methods on a comprehensive set of scenarios, we generate random instances based on a set of (i) power network sizes and structures, (ii) fault probabilities $p$, (iii) repair time durations $s$, and (iv) sizes and (v) shapes of the service region. We present these attributes in Table \ref{Tab:Param}. We generate distribution networks over rectangular and circular regions, using the degree-constrained minimum spanning tree ($d$-MST) algorithm \citep{boldon1996minimum}. After creating the initial set of networks, we generate modifications to control the depth of the network. 
These modifications are created by re-configuring the precedence relationships for a certain number of randomly selected nodes and connecting them to the parent of their parent node. 
We present the details of our problem instance generation method in  Appendix \ref{Ap:InstanceGeneration}.

\begin{table}[h]
\centering
\caption{Parameter values and attributes used in random instance generation.}
 \scalefont{0.8}
\begin{tabular}{@{} ll @{}}   \toprule
\emph{Parameter/Attribute} & \emph{Values}   \\\midrule
Number of nodes $n$ & 20, 30, 40, 50\\
Fault probability $p$    & 0.25, 0.50, 0.75, 0.90  \\ 
Repair duration $s$ & 0, 1.5, 3\\
Depth of the power network $d$ & Modified by re-configuring  $n/2$ or $n/4$ random nodes\\
Shape of the service region & Square, Rectangle, Circle\\
 \bottomrule
\end{tabular}
\label{Tab:Param}
\end{table}

All algorithms are implemented in the Python programming language, and experiments are performed on a high-performance computation platform with an Intel Xeon E5-2680 v4 (Broadwell), 2.40GHz, 14-core processor. We ran batch jobs for all of our experiments with the following configuration: 1 node, 8 cores per node, and 16 GB of RAM per node in a Linux environment.

For each instance, we use our solution methods to compute the lookup tables that contain the approximate values for the states in the Markov decision process. These values are used to determine the actual route for a given disruption scenario. 
For each instance, we generate 1000 random disruption realizations and determine repair routes using our solution methods ($sNRR$, $SA^1$, $SA^2$ and $SA^3$) and the benchmarks ($PS$ and $NN$).
Using these 1000 realizations and the corresponding routes, we compute the average objective function value and a 95\% confidence interval (CI) for the percentage solution quality gap compared to the overall best-performing solution method over the 1000 realizations.
Within this section, we present results for distribution networks with 20, 30 and 50 nodes, each with two different depths of the network. For brevity, we present additional results for instances with 40 nodes and different shapes of regions in Appendix \ref{App:ExtendedResults} as these results are consistent with those presented in the main document.

Recall that our $sNRR$ solution method is based on an approximate dynamic programming approach, and converges to an optimal solution \citep{powell2007approximate}. This approach is strengthened by structural results to speed up computations. Our $SA^1$, $SA^2$ and $SA^3$ algorithms apply different levels of state aggregation. These algorithms are designed to solve larger problem instances in shorter time but do not guarantee an optimal solution.

\begin{figure}[H]
\centering
\subfloat[Original network with depth $d=11$.]{\label{fig20_1}{\includegraphics[width=0.4\textwidth]{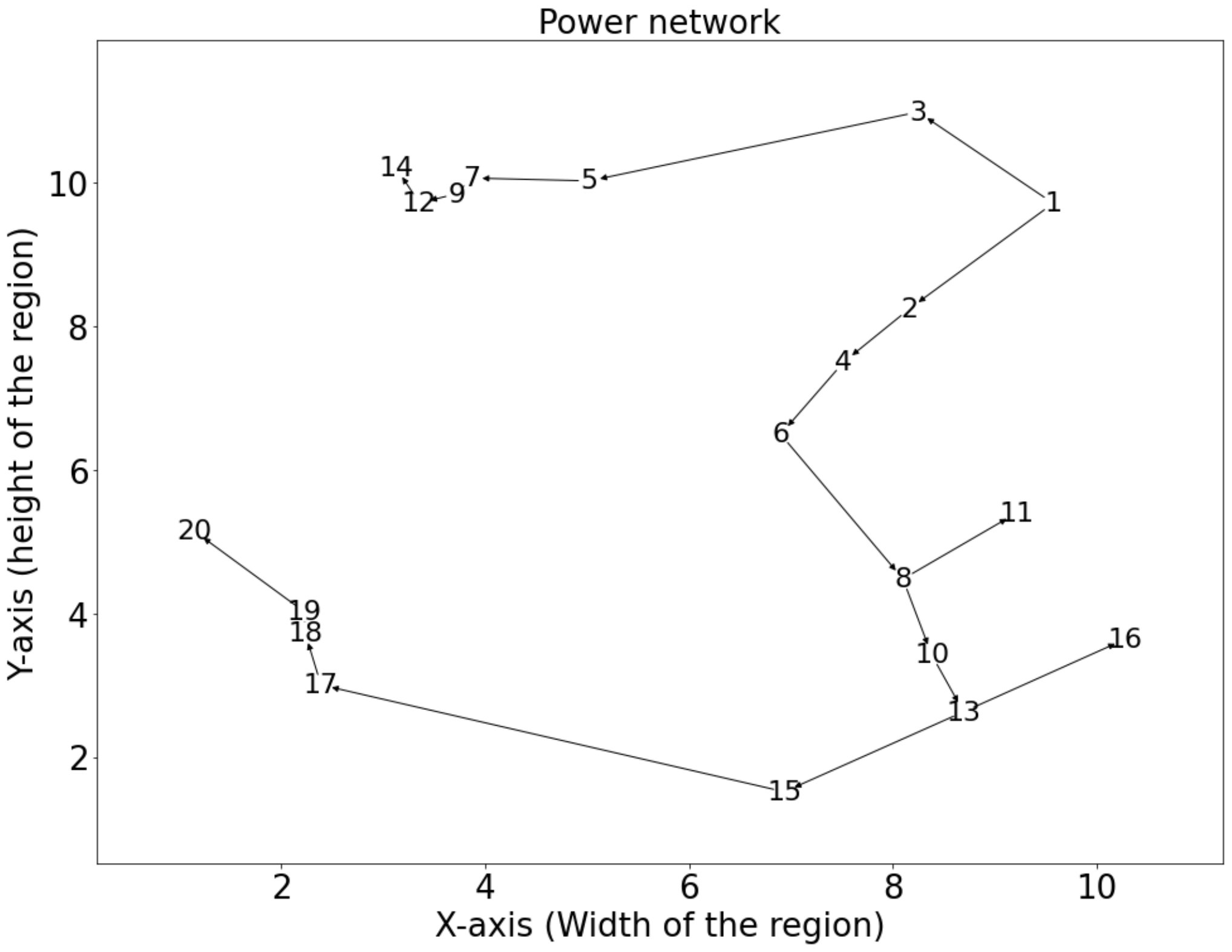}}}
\hspace{2cm}
\subfloat[Modified network with depth $d=6$.]{\label{fig60_2}{\includegraphics[width=0.40\textwidth]{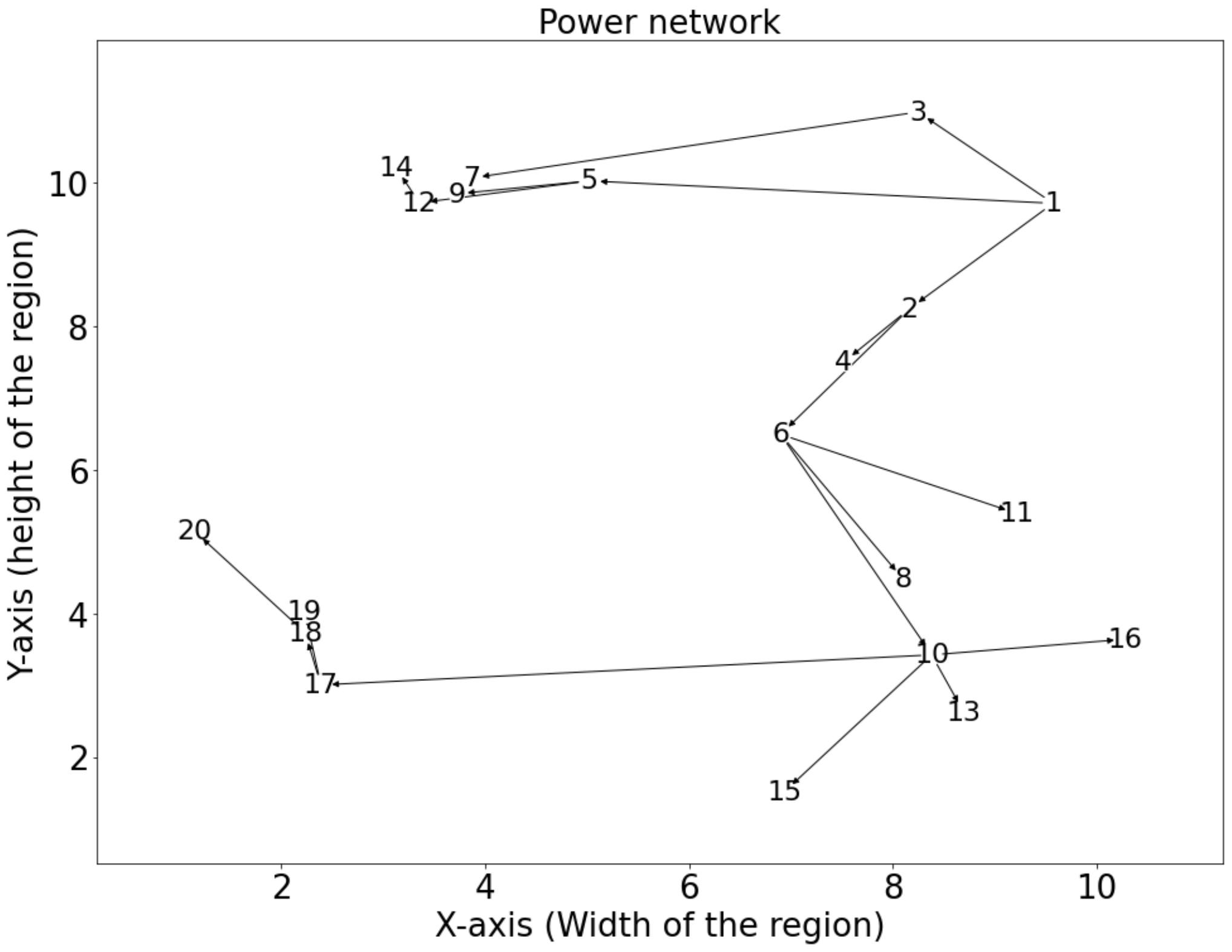}}}
\caption{Node dispersions of power distribution networks with 20 nodes on a $10 \times 10$ square region.}
\label{Fig:20_Square}
\end{figure}

\begin{figure}[h]
\centering
\subfloat[Original network with depth $d=15$.]{\label{fig30_1}{\includegraphics[width=0.40\textwidth]{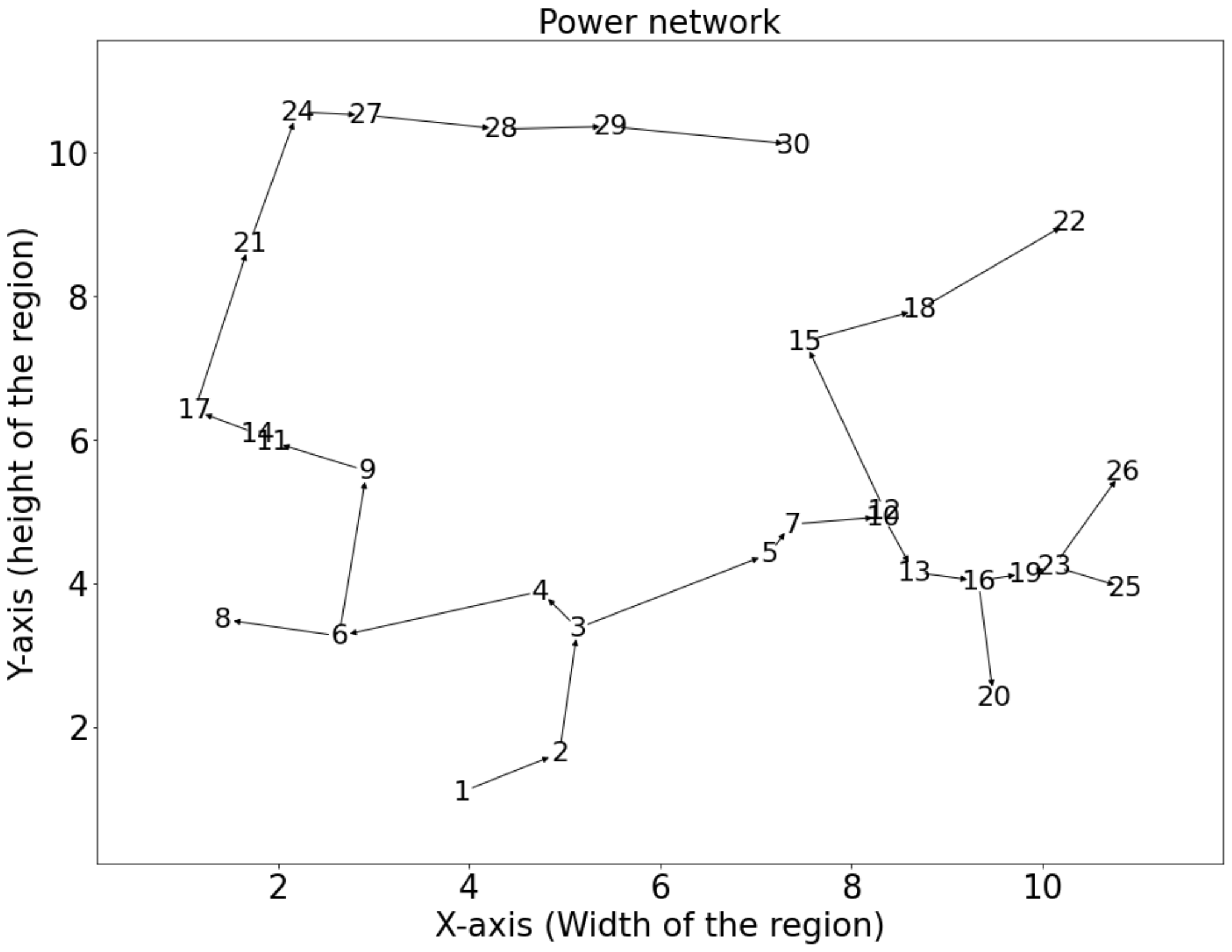}}}
\hspace{1cm}
\subfloat[Modified network with depth $d=11$.]{\label{fig30_4}{\includegraphics[width=0.40\textwidth]{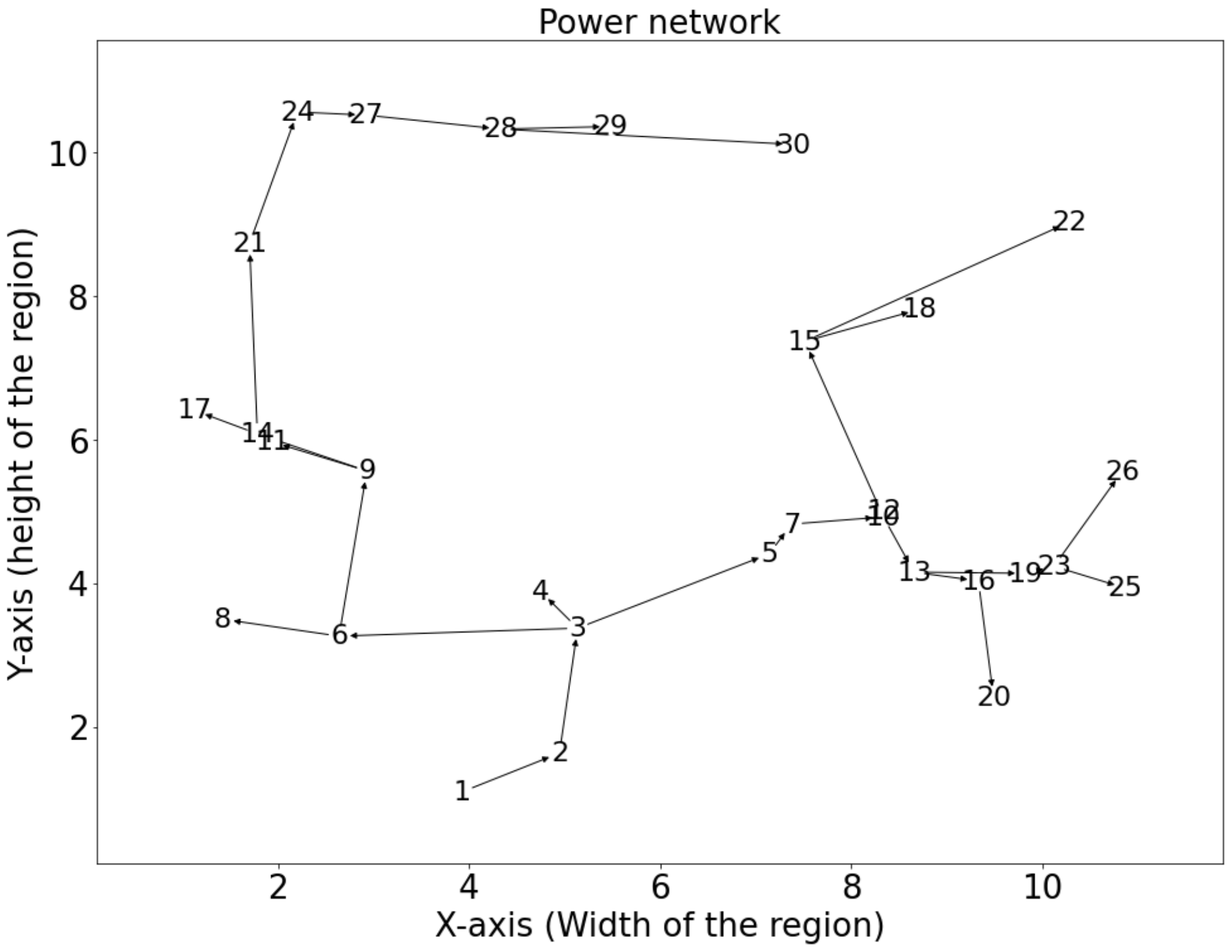}}}
\caption{Node dispersions of power distribution networks with 30 nodes on a 10$\times$10 square region.}
\label{Fig:30_Square}
\end{figure}

\begin{figure}[h]
\centering
\subfloat[Original network with depth $d=15$.]{\label{fig50_1}{\includegraphics[width=0.40\textwidth]{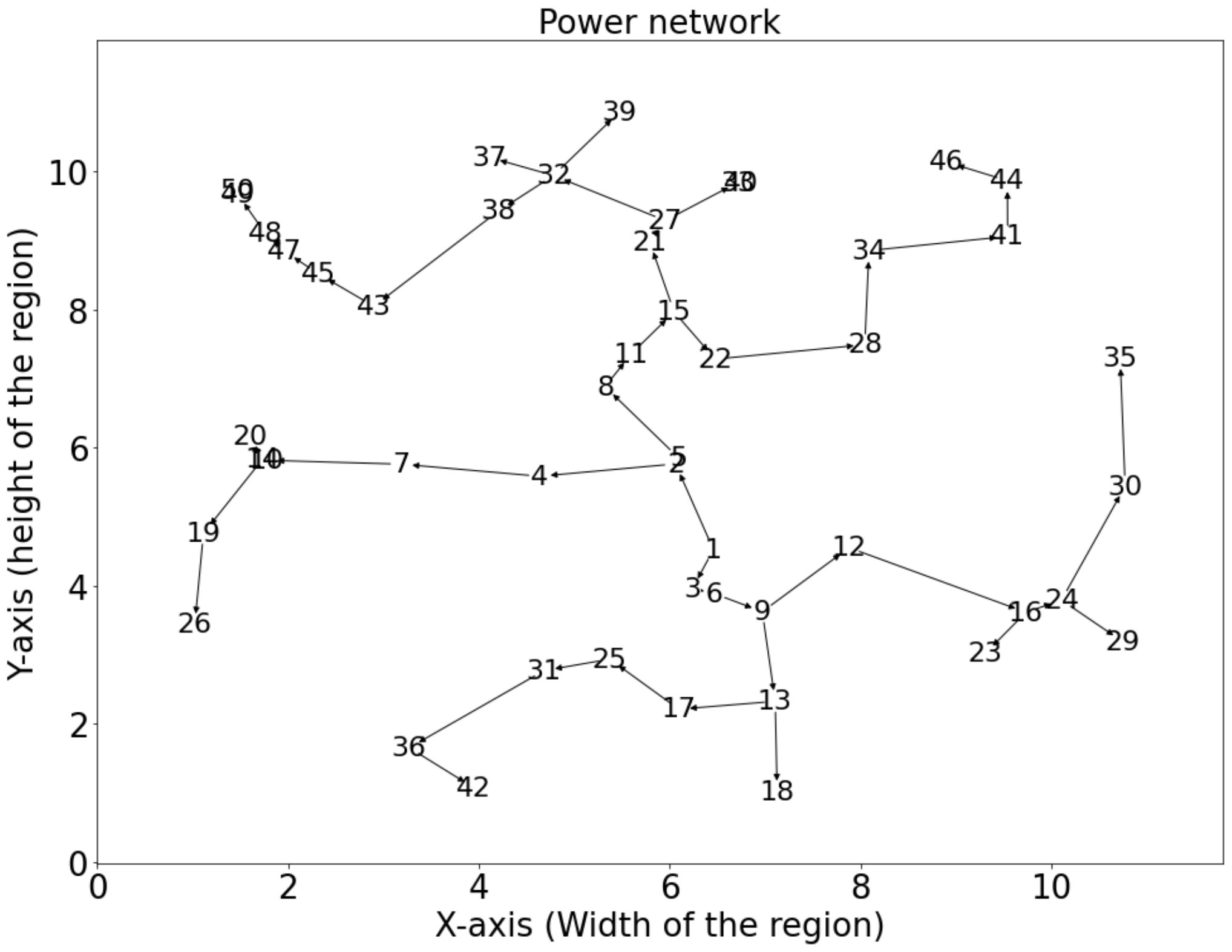}}}\hspace{1cm}
\subfloat[Modified network with depth $d=10$.]{\label{fig50_2}{\includegraphics[width=0.40\textwidth]{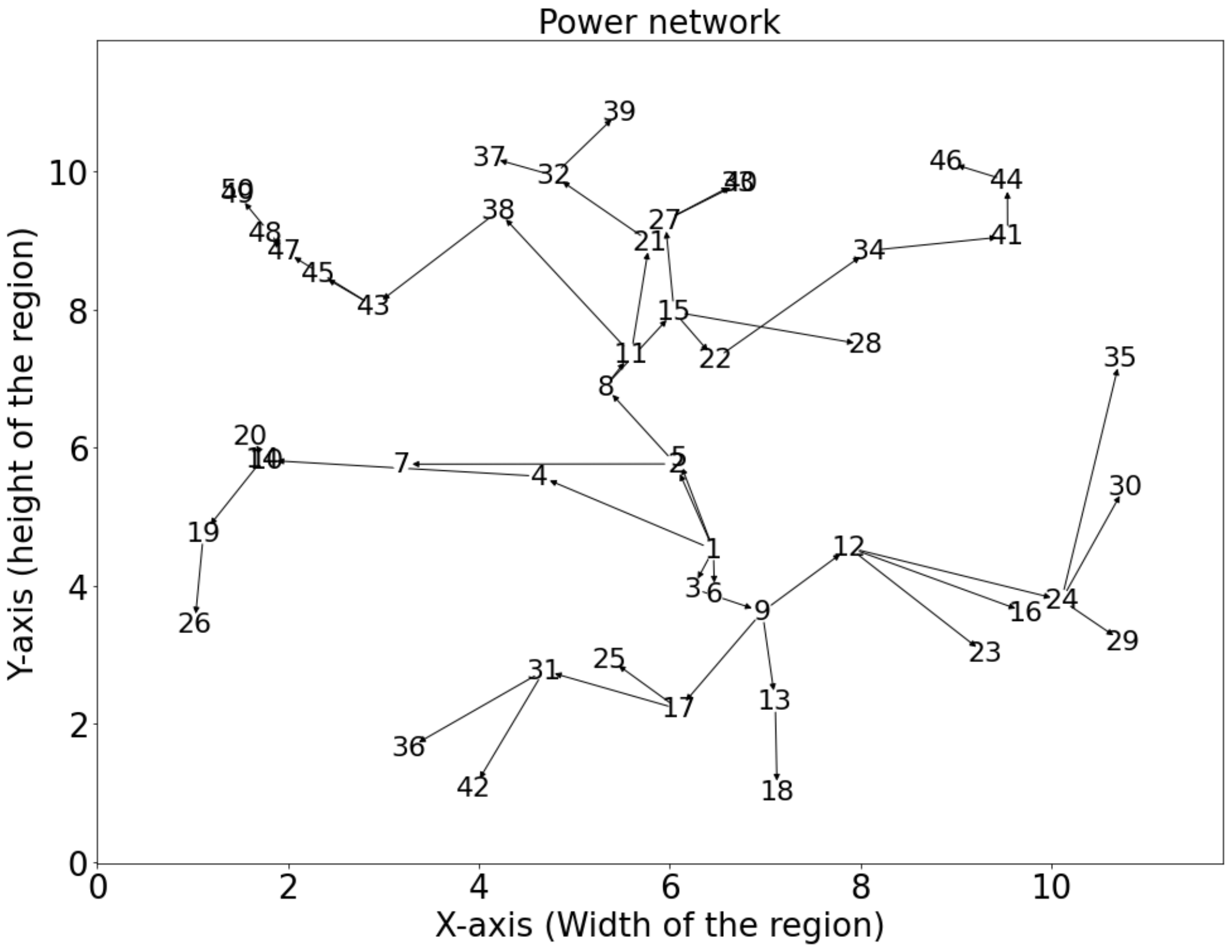}}}
\caption{Node dispersions of power distribution networks with 50 nodes on a 10$\times$10 square region.}
\label{Fig:50_Square}
\end{figure}

Figures \ref{Fig:20_Square}, \ref{Fig:30_Square}, and \ref{Fig:50_Square} illustrate the network layouts on a $10 \times 10$ square region for power distribution networks containing 20, 30, and 50 nodes, respectively.
Tables \ref{Tab:20_sq}, \ref{Tab:30_sq}, and \ref{Tab:50_sq} present the average objective function values over 1000 realizations, as well as the percentage gap between the average solution value for each method and the best-performing method for each set of problem instances for the same distribution networks. 
In the rest of this section, when we comment on the quality of a solution method, we will refer to the average performance over 1000 disruption scenario realizations for a problem instance.

\textbf{Computational performance: } Overall, our $sNRR$ algorithm finds the best solution in all instances, followed by $SA^1$, $SA^2$ and $SA^3$. For example, in Table \ref{Tab:20_sq}, we see that $SA^1$ performs the same as $sNRR$ in 14 of 32 instance sets, deviates by 1\% in 9 instance sets, deviates by 3\% on 7 instance sets, and has a maximum deviation of 7\%. $SA^2$ finds solutions within 3\% of our base $sNRR$ method in 22 of 32 instance sets, and $SA^3$, which uses the highest aggregation level, deviates from $sNRR$ by at most 4\% in 25 of 32 instances. On the other hand, the $PS$ and $NN$ benchmarks perform significantly worse compared to our solution methods. The largest confidence intervals for their percentage gaps compared to $sNRR$ are (29\%, 76\%) and (45\%, 153\%), for $PS$ and $NN$, respectively, for 20 nodes.

The improved solution quality with our $sNRR$, $SA^1$, $SA^2$ and $SA^3$ methods, however, comes with some computational cost.
We present the training times until the stopping criterion is met during approximation and the number of states generated for the MDP process for $n=20$ in Table \ref{Tab:20_sq_ct}. Although the computation time for our $sNRR$ approach is long, we note that this includes the computation time required for generating the look-up tables in the disruption preparation phase. Once these computations are performed, they can be implemented quickly in the response phase to determine the corresponding repair route. Moreover, our aggregate methods ($SA^1$, $SA^2$, and $SA^3$) offer significant reductions both in training time and size of the MDP model. For example, as seen in Table \ref{Tab:20_sq_ct}, for service duration $s=3$, node fault probability $p=0.25$, and network depth $d=11$, $SA^1$, $SA^2$ and $SA^3$ reduce the model training duration by 79\%,  85\% and 99\%, respectively. They also reduce the number of states generated by more than 96\% compared to $sNRR$. Similar observations regarding the relative performance of all algorithms hold across all computations.

\textbf{Fault probability ($p$): } When the repair crew does not have full information regarding the status of a node, i.e., when the node has not yet been visited and has successor nodes without service, the value of $p$ determines the probability of a node fault. Hence, smaller $p$ values imply that it is less likely that a node is faulty, discouraging the repair team from making visits before obtaining more information. In all of the tables presented below, our $sNRR$ algorithm generates the best solutions for all $p$ values. 
We do not observe a significant relationship between the $p$
 value and the performance of our aggregation methods. $SA^1$, $SA^2$ and $SA^3$ find solutions within 4\% of $sNRR$ in a majority of instances for $n=20$. When it comes to the benchmarks, we observe a different behavior. As $p$ increases, the performance of $PS$ decreases while $NN$ performs better. This occurs because, when $p$ is large, the risk of unnecessary trips is reduced.
Hence, the $NN$ heuristic, which prioritizes distance over the priority of nodes, finds better solutions. When $p$ is small, there is a higher risk of unnecessary trips, which can be entirely avoided by the $PS$ heuristic; hence it finds relatively better solutions. However, overall, these myopic benchmarks perform significantly worse than our solution methods, as can be seen by the average percentage gaps in Tables \ref{Tab:20_sq}, \ref{Tab:30_sq} and \ref{Tab:50_sq}.

\begin{table}[h]
 \scalefont{0.8}
\centering
\caption{Comparison of average (Avg) objective value and CI for \% gap compared to $sNRR$ solutions for 20 nodes on a 10$\times$10 square region for distribution networks with different depths.}
\label{Tab:20_sq}
\resizebox{\textwidth}{!}{%
\begin{tabular}{l|l|c c c c |c c ||c c c c|c c}\hline\hline
\multicolumn{2}{l|}{}                            & \multicolumn{6}{c||}{\textbf{Network depth $d = 11$}}                                                                                              
& \multicolumn{6}{c}{\textbf{Network depth $d = 6$}}                                                        \\ 
    \hline 
    \multicolumn{2}{l|}{} & \multicolumn{6}{c||}{\textbf{Avg Total Service Disruption Times}} &   
    \multicolumn{6}{c}{\textbf{Avg Total Service Disruption Times}}  \\ \hline
    $s$ & $p$ & \multicolumn{4}{c|}{\textit{Our Methods}} &  \multicolumn{2}{c||}{\textit{Benchmarks}}   & 
     \multicolumn{4}{c|}{\textit{Our Methods}}  &  \multicolumn{2}{c}{\textit{Benchmarks}}   \\   
    & & $sNRR$  & $SA^1$ &$SA^2$ &$SA^3$ & $PS$ & $NN$ & 
    $sNRR$  & $SA^1$ &$SA^2$ &$SA^3$ & $PS$ & $NN$ \\ \hline 
    \multirow{4}{*}{0} & 0.25 & 149 & 149 & 155 & 155 & 177 & 594   
    & 110 & 110 & 115 & 115 & 125 & 185  \\ 
    & 0.50 & 246 & 246 & 258 & 263 & 387 & 659  
    & 200 & 200 & 210 & 212 & 315 & 295 \\ 
    & 0.75 & 327 & 328 & 338 & 352 & 661 & 701 
    & 275 & 278 & 288 & 295 & 582 & 475  \\ 
    & 0.90 & 360 & 367 & 366 & 396 & 854 & 720  
    & 305 & 315 & 325 & 350 & 850 & 560 \\ \hline 
    \multicolumn{2}{l|}{\multirow{1}{*}{\textbf{CI for Avg \% Gap}}}   
    & -   & (1, 2)  & (2, 3) & (7, 8)  & (73, 76) 
    & \multicolumn{1}{l||}{(60, 63)} & -   & (1, 2)  & (3, 4) & (7, 8)  & (71, 74) &  \multicolumn{1}{l}{(46, 48)} \\ \hline
    \multirow{4}{*}{1.5} & 0.25 & 230 & 230 & 236 & 232 & 252 & 746  
    & 195 & 195 & 200 & 210 & 220 & 290 \\ 
    & 0.50 & 406 & 407 & 418 & 414 & 531 & 949  
    & 350 & 352 & 360 & 370 & 490 & 420 \\ 
    & 0.75 & 563 & 573 & 576 & 582 & 886 & 1131  
    & 525 & 535 & 540 & 555 & 940 & 680 \\ 
    & 0.90 & 671 & 669 & 677 & 681 & 1132 & 1233  
    & 575 & 580 & 590 & 605 & 1200 & 900  \\ \hline 
    \multicolumn{2}{l|}{\multirow{1}{*}{\textbf{CI for Avg \% Gap }}} & -  & (0, 0)  & (3, 3) & (2, 2)  & (52, 54) & \multicolumn{1}{r||}{(150, 153)} 
    & -   & (1, 2)  & (3, 4) & (4, 5)  & (48, 50) & \multicolumn{1}{r}{(53, 56)} \\ \hline
    \multirow{4}{*}{3} & 0.25 & 307 & 307 & 315 & 311 & 327 & 898  
    & 260 & 265 & 270 & 275 & 295 & 350  \\ 
    & 0.50 & 553 & 565 & 574 & 564 & 675 & 1239  
    & 500 & 505 & 510 & 520 & 710 & 620 \\ 
    & 0.75 & 802 & 820 & 828 & 815 & 1110 & 1560  
    & 745 & 755 & 765 & 775 & 1210 & 850 \\ 
    & 0.90 & 944 & 962 & 966 & 956 & 1409 & 1746  
    & 820 & 830 & 850 & 865 & 1590 & 1050  \\ \hline 
    \multicolumn{2}{l|}{\multirow{1}{*}{\textbf{CI for Avg \% Gap }}} & -  & (2, 2)  & (2, 3) & (1, 2)  & (29, 31) & \multicolumn{1}{r||}{(132, 135)}  
    & -   & (1, 1)  & (3, 4) & (3, 4)  & (29, 31) & \multicolumn{1}{r}{(45, 46)} \\ 
    \hline \hline 
\end{tabular}
}
\end{table}

\begin{table}[h]
\centering
 \scalefont{0.8}
\caption{Training duration for value approximation ($t$ in hours) and the number of states (in thousands) in the MDP models for instances with a 20 nodes power network on a 10$\times$10 square region.}
\label{Tab:20_sq_ct}
\resizebox{\textwidth}{!}{%
\begin{tabular}{l|l|c c c c|c c c c|c c c c|c c c c}
\hline\hline
\multicolumn{2}{l|}{}                            & \multicolumn{8}{c|}{\textbf{Network depth $d = 11$}}                                                       
& \multicolumn{8}{c}{\textbf{Network depth $d = 6$}}                                                      \\ \hline 
\multicolumn{2}{l|}{} & \multicolumn{4}{c|}{\textit{t(hr)}} &  \multicolumn{4}{c|}{Number of States} &  
\multicolumn{4}{c|}{\textit{t(hr)}}  & \multicolumn{4}{c}{Number of States} \\ \hline
\multicolumn{1}{l|}{$s$} & $p$ & $sNRR$  & $SA^1$ &$SA^2$ & \multicolumn{1}{l|}{$SA^3$} & $sNRR$ & $SA^1$  &$SA^2$ & $SA^3$ & 

$sNRR$  & $SA^1$ &$SA^2$ &$SA^3$ & $sNRR$ & $SA^1$  &$SA^2$ & $SA^3$ \\ \hline
\multirow{4}{*}{0} & 0.25 & 16.90 & 4.62 & 2.68 & 1.37 & 143.1 & 35.1 & 5.8 & 1.0 

& 14.50 & 2.00 & 2.10 & 1.30 & 284.5 & 65.0 & 9.0 & 1.4 \\
& 0.50 & 31.78 & 4.43 & 3.55 & 2.28 & 515.3 & 41.7 & 6.8 & 1.0 

& 28.50 & 3.20 & 3.80 & 2.70 & 890.5 & 70.0 & 9.2 & 1.4 \\
& 0.75 & 43.82 & 3.95 & 2.81 & 2.77 & 1,415.2 & 44.7 & 7.1 & 1.0  
& 42.00 & 4.00 & 3.50 & 2.60 & 1,680.0 & 74.5 & 9.5 & 1.5 \\
& 0.90 & 22.14 & 3.64 & 2.90 & 2.71 & 2,066.4 & 45.7 & 7.2 & 1.0  
& 20.00 & 3.60 & 3.20 & 2.50 & 2,115.0 & 75.0 & 9.5 & 1.5 \\ \hline
\multirow{4}{*}{1.5} & 0.25 & 15.73 & 3.43 & 2.70 & 1.37 & 563.8 & 45.0 & 7.0 & 1.1 
& 16.50 & 2.30 & 2.15 & 1.35 & 856.0 & 75.2 & 9.5 & 1.4 \\
& 0.50 & 33.72 & 8.98 & 5.17 & 2.60 & 1,713.6 & 47.5 & 7.3 & 1.1 
& 31.50 & 4.50 & 4.10 & 2.70 & 2,490.0 & 75.5 & 9.6 & 1.4 \\
& 0.75 & 47.66 & 6.02 & 5.99 & 2.23 & 4,529.2 & 49.2 & 7.5 & 1.1  
& 45.00 & 4.00 & 3.80 & 2.60 & 6,240.0 & 75.0 & 9.7 & 1.4 \\
& 0.90 & 54.55 & 3.15 & 3.01 & 2.57 & 6,490.5 & 49.7 & 7.5 & 1.0  
& 50.00 & 3.50 & 3.30 & 2.50 & 7,890.0 & 75.5 & 9.7 & 1.5 \\ \hline
\multirow{4}{*}{3} & 0.25 & 16.54 & 3.45 & 2.41 & 1.37 & 1,208.7 & 49.3 & 7.6 & 1.1 
& 17.00 & 2.30 & 2.15 & 1.35 & 1,750.0 & 75.0 & 9.5 & 1.4 \\
& 0.50 & 34.70 & 5.98 & 4.95 & 2.66 & 3,559.5 & 50.6 & 7.6 & 1.1 
& 34.00 & 4.50 & 4.20 & 2.70 & 4,850.0 & 75.5 & 9.6 & 1.4 \\
& 0.75 & 50.73 & 6.08 & 4.42 & 3.79 & 7,921.9 & 51.3 & 7.7 & 1.1 
& 48.00 & 6.50 & 5.90 & 3.80 & 9,440.0 & 75.5 & 9.7 & 1.4 \\
& 0.90 & 57.07 & 6.45 & 4.24 & 2.62 & 9,678.1 & 51.5 & 7.7 & 1.1 
& 52.00 & 5.90 & 5.40 & 3.70 & 10,650.0 & 75.5 & 9.7 & 1.5 \\ 
\hline \hline
\end{tabular}
}
\end{table}

\textbf{Repair duration ($s$): } Our problem shares similarities with scheduling under precedence constraints and the TRP. When the repair duration $s$ is small, repair crew routes become similar to TRP routes, whereas for larger $s$, scheduling decisions become more dominant, rendering the travel durations less important.  Hence, when $s$ is large enough compared to the average travel time between nodes, the repair crew routing problem becomes easier to solve, and the optimal solution follows the order of priority.
This is evident by the improving performance of the $PS$ benchmark heuristic for large $s$ values as seen in Tables \ref{Tab:20_sq}, \ref{Tab:30_sq} and \ref{Tab:50_sq}. Similarly, the gap between the $NN$ and the best available solution decreases. However, this is not because the $NN$ heuristic performs better, but because the additional travel time caused by the $NN$ algorithm has a lower impact on the entire objective function value when the service duration is large. On the other hand, our solution methods outperform the benchmarks significantly. 
Overall, the performance of our solution methods is not affected by the $s$ value. 
However, we observe a slight decrease in the quality of routes determined by the $SA^3$ method when $s$ is small and $p$ is large. When $s$ is small, the repair times are mainly determined by the time spent during routing and less by the time spent while fixing faults. Therefore, when $s$ is smaller, routing decisions become more important. On the other hand, for larger $s$ values, it is the opposite, and the precedence relationships on the distribution network become more dominant when determining the travel sequence. Therefore, for small $s$ values, the objective function value depends more on the quality of the route rather than the precedence sequence.
Moreover, when the probability $p$ is high, there is a larger amount of information available, which is aggregated more heavily by the $SA^3$ method, which can lead to lower-quality routing decisions.

\textbf{Network depth ($d$): } As the depth of the infrastructure network increases, the average number of immediate successor nodes decreases, and the revelation of information on node status (i.e., faulty or not) becomes slower. Therefore, for instances with larger distribution network depth and the same node locations, we observe that the total service disruption durations increase in general as we do not obtain much actionable information. 
This trend can be observed in Tables \ref{Tab:20_sq}, \ref{Tab:30_sq}, and \ref{Tab:50_sq}. 
Consequently, a reduced network depth leads to more effective repair strategies and lower associated costs. However, as seen in Table \ref{Tab:50_sq}, for example, when the fault probability is high (causing slower information revelation)  and service times are short relative to travel distances, an increase in network depth does not necessarily result in higher average objective values for the $PS$ and $NN$ heuristics.

\begin{table}[H]
\centering
 \scalefont{0.8}
\caption{Comparison of average (Avg) objective values and \% CI's for gap compared to $SA^2$ solutions for 30 nodes on a 10$\times$10 square region for different depths.}
\label{Tab:30_sq}
\resizebox{0.8\textwidth}{!}{\begin{tabular}{l|l|c c|c c|c c |c c }
\hline \hline
\multicolumn{2}{l|}{}                            & \multicolumn{4}{c|}{Network depth $d = 15$}                                                      & \multicolumn{4}{c}{Network depth $d =11$}                                                        \\ \hline\multicolumn{2}{l|}{} & \multicolumn{4}{c|}{\makecell{\textbf{Avg Total Service } \\ \textbf{Disruption Times}}} & \multicolumn{4}{c}{\makecell{\textbf{Avg Total Service } \\ \textbf{Disruption Times}}} \\ \hline
$s$ & $p$ & \multicolumn{2}{c|}{\textit{Our Methods}} &  \multicolumn{2}{c|}{\textit{Benchmarks}}   & 
     \multicolumn{2}{c|}{\textit{Our Methods}}  &  \multicolumn{2}{c}{\textit{Benchmarks}} \\
                 &  & $SA^2$  & $SA^3$  &$PS$   & \multicolumn{1}{l|}{NN}   & $SA^2$   & $SA^3$  & $PS$    & \multicolumn{1}{l}{NN}\\ \hline
\multicolumn{1}{l|}{\multirow{4}{*}{0}}   & 0.25 & 383  & 384  & 417  & \multicolumn{1}{l|}{1708}  & 345  & 346  & 387  & \multicolumn{1}{l}{1708}  \\
\multicolumn{1}{l|}{}                     & 0.5  & 525  & 529  & 654  & \multicolumn{1}{l|}{1708} & 502  & 507  & 686  & \multicolumn{1}{l}{1708}  \\
\multicolumn{1}{l|}{}                     & 0.75 & 634  & 647  & 912  & \multicolumn{1}{l|}{1708}  & 628  & 642  & 1014 & \multicolumn{1}{l}{1708}  \\
\multicolumn{1}{l|}{}                     & 0.9  & 689  & 713  & 1100 & \multicolumn{1}{l|}{1708} & 695  & 715  & 1241 & \multicolumn{1}{l}{1708} \\ \hline
\multicolumn{2}{l|}{\multirow{1}{*}{\textbf{CI for Avg \% Gap }}} & -  & (1, 2)  & (33, 35)  & \multicolumn{1}{l|}{(225, 230)}  & -  & (0, 1)  & (47, 49)  & \multicolumn{1}{l}{(250, 253)}  \\ \hline
\multicolumn{1}{l|}{\multirow{4}{*}{1.5}} & 0.25 & 563  & 562  & 584  & \multicolumn{1}{l|}{2073}  & 497  & 496  & 523  & \multicolumn{1}{l}{2073} \\
\multicolumn{1}{l|}{}                     & 0.5  & 874  & 872  & 977  & \multicolumn{1}{l|}{2407}& 828  & 827  & 965  & \multicolumn{1}{l}{2407}  \\
\multicolumn{1}{l|}{}                     & 0.75 & 1165 & 1167 & 1413 & \multicolumn{1}{l|}{2731} & 1147 & 1150 & 1481 & \multicolumn{1}{l}{2731}  \\
\multicolumn{1}{l|}{}                     & 0.9  & 1333 & 1338 & 1717 & \multicolumn{1}{l|}{2927}  & 1326 & 1338 & 1840 & \multicolumn{1}{l}{2927} \\ \hline
\multicolumn{2}{l|}{\multirow{1}{*}{\textbf{CI for Avg \% Gap  }}} & -  & (0, 0)  & (17, 19)  & \multicolumn{1}{l|}{(180, 183)} & -  & (0, 0)  & (23, 25)  & \multicolumn{1}{l}{(200, 202)}  \\ \hline
\multicolumn{1}{l|}{\multirow{4}{*}{3}}   & 0.25 & 735  & 733  & 750  & \multicolumn{1}{l|}{2439}  & 636  & 635  & 658  & \multicolumn{1}{l}{2439}  \\
\multicolumn{1}{l|}{}                     & 0.5  & 1213 & 1210 & 1299 & \multicolumn{1}{l|}{3107} & 1133 & 1132 & 1244 & \multicolumn{1}{l}{3107}  \\
\multicolumn{1}{l|}{}                     & 0.75 & 1683 & 1685 & 1915 & \multicolumn{1}{l|}{3753}  & 1656 & 1656 & 1947 & \multicolumn{1}{l}{3753}\\
\multicolumn{1}{l|}{}                     & 0.9  & 1965 & 1963 & 2333 & \multicolumn{1}{l|}{4146}  & 1961 & 1960 & 2440 & \multicolumn{1}{l}{4146}\\ \hline
\multicolumn{2}{l|}{\multirow{1}{*}{\textbf{CI for Avg \% Gap  }}} & -  & (0, 0)  & (9, 12)  & \multicolumn{1}{l|}{(160, 162)}  & -  & (0, 0)  & (16, 18)  & \multicolumn{1}{l}{(175, 178)} \\ 
\hline \hline
\end{tabular}}
\end{table}

\begin{table}[H]
\centering
 \scalefont{0.8}
\caption{Training duration for value approximation ($t$ in hours) and the number of states (in thousands) in the MDP models for instances with a 30 nodes power network on a 10$\times$10 square region.}
\label{30_sq_CT}
\resizebox{0.7\textwidth}{!}{
\begin{tabular}{ll|cc|cc|cc|cc}
\hline \hline
\multicolumn{2}{l|}{\multirow{2}{*}{}}           & \multicolumn{4}{c|}{\textbf{Network depth} $d=15$}                                                       & \multicolumn{4}{c}{\textbf{Network depth} $d=11$}                                                         \\ \cline{3-6} \cline{7-10}
\multicolumn{2}{l|}{}                            & \multicolumn{2}{c|}{\textit{t(hr)}}     & \multicolumn{2}{c|}{Number of States} & \multicolumn{2}{c|}{\textit{t(hr)}}     & \multicolumn{2}{c}{Number of States} \\ \hline
\multicolumn{1}{l|}{s}                    & $p$    & $SA^2$   & $SA^3$    & $SA^2 $          & $SA^3 $          & $SA^2$   & $SA^3$   & $SA^2 $   & $SA^3 $            \\ \hline
\multicolumn{1}{l|}{\multirow{4}{*}{0}}   & 0.25 & 6.17  & 5.13  &  23.1       & 2.4          & 10.43 & 6.18   & 32.3      & 3.9      \\
\multicolumn{1}{l|}{}                     & 0.5  & 8.83  & 8.11  & 25.1     & 2.4        & 12.85 & 7.60  & 37.8      & 4.5           \\
\multicolumn{1}{l|}{}                     & 0.75 & 11.07 & 9.99   & 26.0       & 2.5          & 14.12 & 10.20  & 42.7      & 4.8       \\
\multicolumn{1}{l|}{}                     & 0.9  & 13.18 & 11.14& 26.2       & 2.5     & 15.21 & 10.59 & 43.4      & 4.9          \\ \hline
\multicolumn{1}{l|}{\multirow{4}{*}{1.5}} & 0.25 & 18.30 & 15.92 & 26.2       & 2.5          & 19.30 & 15.07& 37.8      & 4.2         \\
\multicolumn{1}{l|}{}                     & 0.5  & 21.20 & 20.62& 26.7       & 2.5           & 34.45 & 26.42  & 42.8      & 4.8          \\
\multicolumn{1}{l|}{}                     & 0.75 & 15.26 & 12.13  & 27.0       & 2.5         & 39.09 & 31.86  & 44.4      & 4.9          \\
\multicolumn{1}{l|}{}                     & 0.9  & 18.56 & 11.17 & 27.0       & 2.5         & 21.79 & 15.65  & 44.4      & 4.9           \\ \hline
\multicolumn{1}{l|}{\multirow{4}{*}{3}}   & 0.25 & 21.14 & 16.18  & 27.7       & 2.5           & 20.66 & 11.69  & 38.6      & 4.2          \\
\multicolumn{1}{l|}{}                     & 0.5  & 43.42 & 30.38 & 27.4       & 2.5           & 47.56 & 35.78& 44.6      & 4.9           \\
\multicolumn{1}{l|}{}                     & 0.75 & 21.61 & 18.66 & 27.6       & 2.5         & 48.56 & 42.16 & 45.0      & 4.9         \\
\multicolumn{1}{l|}{}                     & 0.9  & 13.94 & 11.17  & 27.6       & 2.5           & 22.07 & 15.33  & 44.5      & 4.9  \\       
\hline\hline
\end{tabular}}

\end{table}

\textbf{Network size: }
As the network size increases, state aggregation methods, particularly $SA^3$, show significant reductions in both training duration and the size of the model compared to $sNRR$. In smaller networks, as indicated by Table \ref{Tab:20_sq_ct}, $SA^1$ and $SA^2$ perform similarly. Although $SA^1$ and $SA^2$ slightly compromise solution quality compared to $sNRR$, they reduce the model training duration on average by 80\%.  
For networks with 30 nodes, $SA^3$ offers a significant computational advantage over $SA^2$, with an average of a 50\% reduction in training time and a 91\% decrease in the number of states. For networks with 40 and 50 nodes (Tables \ref{Tab:40_sq} and \ref{Tab:50_sq}), only $SA^3$ is feasible due to memory constraints. 

\begin{table}[H]
\centering
 \scalefont{0.8}
\caption{Comparison of average (Avg) objective values and \% CI's for gap compared to $SA^3$ solutions for 50 nodes on a 10$\times$10 square region for different depths.}
\label{Tab:50_sq}

\resizebox{0.8\textwidth}{!}{
\begin{tabular}{l|l|r|rr|r|rr}
\hline \hline 
\multicolumn{2}{c|}{} & \multicolumn{3}{c|}{\textbf{Network depth} $d=15$} 
& \multicolumn{3}{c}{\textbf{Network depth} $d=10$} \\ \hline
\multicolumn{2}{l|}{} & \multicolumn{3}{c|}{\makecell{\textbf{Avg Total Service } \\ \textbf{Disruption Times}}} 
& \multicolumn{3}{c}{\makecell{\textbf{Avg Total Service} \\ \textbf{ Disruption Times}}}\\ \hline
      $s$ & $p$ & \multicolumn{1}{c|}{\textit{Our Methods}} &  \multicolumn{2}{c|}{\textit{Benchmarks}}   & 
     \multicolumn{1}{c|}{\textit{Our Methods}}  &  \multicolumn{2}{c}{\textit{Benchmarks}}\\ &  & $SA^3$ & $PS$ & $NN$ 
& $SA^3$ & $PS$ & $NN$ \\ \hline
\multirow{4}{*}{0}   & 0.25 & 600  & 832   & 2530  
& 521  & 692   & 2471  \\
                     & 0.5  & 1034 & 2041  & 2748    
                     & 989  & 1806  & 2779   \\
                     & 0.75 & 1412 & 3560  & 2869   
                     & 1357 & 3315  & 2952    \\
                     & 0.9  & 1659 & 4550  & 2894      
                     & 1535 & 4346  & 3007   \\ \hline
                       \multicolumn{2}{l|}{\multirow{1}{*}{\textbf{CI for Avg \% Gap }}} & -  & (123, 123)  & (178, 179)   
                       & -   & (102, 104)  & (195, 197) \\ \hline
\multirow{4}{*}{1.5} & 0.25 & 979  & 1181  & 3425  
& 846  & 985   & 3351   \\
                     & 0.5  & 1889 & 2832  & 4478   
                     & 1784 & 2508  & 4518   \\
                     & 0.75 & 2844 & 4883  & 5430   
                     & 2759 & 4567  & 5558   \\
                     & 0.9  & 3472 & 6219  & 5933   
                     & 3291 & 5981  & 6117 \\ \hline
                       \multicolumn{2}{l|}{\multirow{1}{*}{\textbf{CI for Avg \% Gap }}} & -  & (58, 59)  & (147, 150) 
                       & -   & (48, 48)  & (163, 165) \\ \hline
\multirow{4}{*}{3}   & 0.25 & 1351 & 1530  & 4321   
& 1163 & 1279  & 4230   \\
                     & 0.5  & 2724 & 3623  & 6207     
                     & 2565 & 3209  & 6258   \\
                     & 0.75 & 4196 & 6206  & 7992  
                     & 4140 & 5819  & 8164   \\
                     & 0.9  & 5171 & 7888  & 8972   
                     & 5153 & 7616  & 9227   \\ \hline
                       \multicolumn{2}{l|}{\multirow{1}{*}{\textbf{CI for Avg \% Gap }}} & -  & (38, 39)  & (133, 134)   
                       & -   & (30, 31)  & (152, 153) \\ \hline

\hline 
\end{tabular}} 
\end{table}

We performed additional computational experiments on instances where nodes are distributed over different area shapes, e.g., rectangular and circular. The observations from these instances are consistent with those we have made for square regions, suggesting that the service region shape does not have a significant impact on the relative performance of solution methods.
More specifically, for problems that are equivalent in network size and precedence structure, as the service region shape changes (which affects the travel times between nodes), our observations are consistent with those when $s$ changes. For example, as the region becomes more elongated, average travel times between nodes increase. Hence the results in such cases are similar to previous observations when $s$ decreases.
Hence, we omit these results for brevity.


\section{Conclusion}
\label{Sec:Conc}

This paper studies repair crew routing for infrastructure network restoration under incomplete information, mainly focusing on power distribution systems. 
Our goal is to determine a crew repair route with minimum total expected service disruption duration. 
We call this problem the Traveling Repairman Network Restoration Problem (TRNRP).  
While the TRNRP is similar to the Traveling Repairman Problem (TRP), it presents additional challenges due to inter-dependencies between infrastructure network nodes and incomplete information on node status (i.e., whether or not a node is faulty).
Considering the dynamic nature of the corresponding routing decisions, we formulate the problem as a Markov Decision Process and use a dynamic programming solution method. 
To address the computational challenges due to the size of the problem, we develop an approximate dynamic programming approach using reinforcement learning. 
In addition, we present structural results to identify a set of suboptimal solutions to improve the model training process. Moreover, we develop new state aggregation methods leveraging information on network inter-dependencies to solve larger problem instances. 
We perform extensive computational experiments to evaluate our solution methods on a rich set of instances based on various infrastructure network topologies and problem parameters. We provide a summary of how our solution methods perform compared to benchmarks under different problem parameter settings. Our methods find significantly better solutions than the benchmarks and require reasonable model training times. In addition, our state aggregation methods reduce the computation time requirements significantly, e.g., by more than 87.5\% in training time and 99\% in required memory.

This work suggests multiple new future research directions.
One key area is adapting the approximate dynamic programming algorithm to problems where the road network can also be disrupted, and the crew has incomplete information on arc status. In addition, making dynamic crew routing decisions for multiple crews represents an open challenge.  Finally, future work may consider network switching decisions that can both dynamically change the network's precedence structure and reveal hidden information during repair operations.

\singlespacing

\bibliographystyle{apalike}

\let\oldthebibliography\thebibliography
\let\endoldthebibliography\endthebibliography
\renewenvironment{thebibliography}[1]{
  \begin{oldthebibliography}{#1}
    \setlength{\itemsep}{0em}
    \setlength{\parskip}{0.5em}
}
{
  \end{oldthebibliography}
}

\bibliography{Power.bib}

\doublespacing

\newpage
\appendix
\section{Algorithms}
\label{Ap:TransProp}
\begin{algorithm}[h]
\caption{Calculating Transition Probabilities}\label{Alg:TransitionPseudo}
\resizebox{0.8\textwidth}{!}{
\begin{minipage}{\textwidth}
 \begin{algorithmic}[1]

\State $\beta_1, \beta_2, \beta_3, p_1, p_2 \gets 1$ 
\State $U^{+}:$ set of nodes with service
\State $U^{-}_1:$  set of nodes without service and faulty with probability 1  
\State    $U^{-}_0:$ set of nodes without service and faulty with probability 0  
\State    $U^{-}_p: $ set of nodes without service and faulty with probability $p$
\State    $a_t:$  is the action 

\If{$a_t \in U^{-}_1$} 
    \State Node $a_t$ moves to the set $U^{+}$ 
    \For {$j \in U^{-}_1 \cup U^{-}_0 \cup U^{-}_p$} 
    \If{$a_t \notin P_{1j}$} \Comment{node $a_t$ is not a predecessor of node $j$}
    \State $j$ stays in the same set
    \State $\beta_1 = \beta_1 \times 1$ \Comment{Updating transition probability factor}
    \ElsIf{$a_t \in P_{1j}$} \Comment{node $a_t$ belongs to the unique path between the source and $j$}
    \If {$(a_t, j) \in A_{\mathcal{I}}$ and $j$ gains power} \Comment{If there is an edge between $a_t$ and $j$}
    \State $\beta_2 = \beta_2 \times (1 - p)$ \Comment{Node $j$ is not faulty with prob $(1-p)$}
    \State Node $j$ moves to the set $U^{+}$ 
    \For {$h \in U^{-}_1 \cup U^{-}_0 \cup U^{-}_p$} 
    \If{$j\preceq h$ and $h$ gains power} \Comment{Checking if node $h$ is a successor of node $j$}
    \State $\beta_3 = \beta_3 \times (1 - p)$ \Comment{Node $h$ is not faulty with prob $(1-p)$}
    \State Node $h$ moves to the set $U^{+}$ 
    \ElsIf {$j\preceq h$ and $h$ doesn't gain power}
    \State $\beta_3 = \beta_3 \times p$ \Comment{Node $h$ is faulty with prob $p$}
    \State Node $h$ moves to the set $U^{-}_1$
    \EndIf
    \EndFor
    \ElsIf{$j$ doesn't gain power} 
    \State $\beta_2 = \beta_2 \times p$ \Comment{Node $j$ is faulty with prob $p$}
    \State Node $j$ moves to the set $U^{-}_1$ 
    \EndIf
    \EndIf
    \EndFor
    \State $p_1 = \beta_1 \times \beta_2 \times \beta_3$ \Comment{Calculating combined probability when the action $a_t \in U^{-}_1$}
    \ElsIf {$ a_t \in U^{-}_p$} 
    \For{$j \in P_{1a_t}$ and $j \in U^{-}_1 \cup U^{-}_p$} 
    \State $p_2 = p_2 \times p$ \Comment{Node $j$ is faulty with prob $p$}
    \EndFor
\EndIf
\State $P(S_{t+1}|S_t, a_t) = p_1 \times p_2$ \Comment{Calculating final transition probability}

\end{algorithmic}
\end{minipage}}
\end{algorithm}

\begin{algorithm}

\caption{: NRR Algorithm}\label{alg:cap1}
\resizebox{0.7\textwidth}{!}{
\begin{minipage}{\textwidth}
\begin{algorithmic}

\State $\delta \gets$ probability values to determine the exploration-exploitation trade-off
\State $L \gets$ Depot
\State $G_{\mathcal{I}} \gets$ Distribution Network 
\State $\|N_{\mathcal{I}}\|\gets$ Number of nodes without power
\State $p \gets$ Probability of a node being faulty 
\State $L \gets$ Current Location 
\State $t$: Decision epoch
\State Post-decision state $(S^{a}_t) \gets  \{a, U^{+}, U^{-}_0, U^{-}_{1}, U^{-}_p \}$
\State Pre-decision state $(S_{t+1}) \gets\{U^{+}, U^{-}_0, U^{-}_{1}, U^{-}_p \}$
\For {$n$ = $1, 2, \dots $} \Comment{outer loop for iterations}
\State $t \gets 0$
\State $S^{a,n}_t=\{\{L\}, \{\},\{1\}, \{\}, \{2, 3, \dots N\}\}$ \Comment{Initial post-decision state} 
\State $\hat{H}^{a,n}_{t}(S^{a,n}_t) = 0$ 
\Comment{Initializing value of initial post-decision state} 
\State $S_{t+1}=\{\{\},\{1\}, \{\}, \{2, 3, \dots N\}\}$ 
\Comment{Initial pre-decision state}
\State $S_F = \{\{1, 2, 3, \dots N\}, \{\}, \{\}, \{\}\}$ \Comment{Final state}

\State retrieve $U^{+}, U^{-}_0, U^{-}_{1}, U^{-}_p$ from pre-decision state $S_{t+1}$
\If {($rand(0, 1) \geq \delta$)}
\While {($S_{t+1} \neq S_F$)} \Comment{inner loop for decision epochs using pure exploitation for $\delta = 0$}
\State $a_t = \argmin \{c(S_t, a_t) + \hat{H}^{a,n-1}_{t}(S^{a}_t)\}$

\State $\hat{h}^n_{t}(S_t^{a, n}) = \{c(S_t, a_t) + \hat{H}^{a,n-1}_{t}(S^{a}_t)\}$
\State $t \gets t+1$
\State $L \gets a_{t}$
\State \{$S^{a, n}_{t}, S_{t+1} \}\gets$ \Call{State Calculation}{$L, G_{\mathcal{I}}, p, U^{+}, U^{-}_0, U^{-}_{1}, U^{-}_p $}
\Comment{This will update $U^{+}, U^{-}_0, U^{-}_{1}, U^{-}_p$ in $S_{t+1}$}
\State $N(S^{a, n}_t) \gets $ Number of times $S^{a, n}_t$ is observed till iteration $n$
\State $\alpha \gets \frac{1}{N(S^{a, n}_t)}$

\State $\hat{H}^{a,n}_{t}(S^{a,n}_t) = (1-\alpha_{n-1})\hat{H}^{a,n-1}_{t}(S^{a,n}_t) + \alpha_{n-1} \hat{h}^n_{t}(S_t^{a, n})$
\EndWhile
\ElsIf{{($rand(0, 1) \leq \delta$)}}
\While {($S_{t+1} \neq S_F$)} \Comment{inner loop for decision epochs using pure exploration for $\delta = 1$}
\State $a_t = $ pick a random node from the set of nodes without power

\State $\hat{h}^n_{t}(S_t^{a, n}) = \{c(S_t, a_t) + \hat{H}^{a,n-1}_{t}(S^{a}_t)\}$
\State $t \gets t+1$
\State $L \gets a_{t}$
\State \{$S^{a, n}_{t}, S_{t+1} \}\gets$ \Call{State Calculation}{$L, G_{\mathcal{I}}, p, U^{+}, U^{-}_0, U^{-}_{1}, U^{-}_p $}
\Comment{This will update $U^{+}, U^{-}_0, U^{-}_{1}, U^{-}_p$ in $S_{t+1}$}
\State $N(S^{a, n}_t) \gets $ Number of times $S^{a, n}_t$ is observed till iteration $n$
\State $\alpha \gets \frac{1}{N(S^{a, n}_t)}$

\State $\hat{H}^{a,n}_{t}(S^{a,n}_t) = (1-\alpha_{n-1})\hat{H}^{a,n-1}_{t}(S^{a,n}_t) + \alpha_{n-1} \hat{h}^n_{t}(S_t^{a, n})$
\EndWhile
\EndIf
\State Stop the \textbf{for} loop once Equation \eqref{stopping} is satisfied.
\EndFor
\Function{State Calculation}{$a_t, G_{\mathcal{I}}, p, U^{+}, U^{-}_1, U^{-}_0, U^{-}_p$} \Comment{function to calculate post-decision state and pre-decision state}
\If {$a_t \in U^{-}_p$}
\State Move node $a_t$ from $U^{-}_p$ to $U^{-}_0$ and update $S^{a}_{t}$
\State No new information to update $S_{t+1}$
\EndIf
\If {$a_t \in U^{-}_1$}
\State Move node $a_t$ from $U^{-}_1$ to $U^{+}$ and update $S^{a}_{t}$
\State {$\vec{\omega} \gets $} stores information after calling  \Call{Simulation}{$a_t, G_{\mathcal{I}}, p $} 
\State use $\vec{\omega}$ to update $S_{t+1}$
\EndIf
\State \Return $S^{a}_{t}, S_{t+1}$
\EndFunction
\Function{Simulation}{$a_t, G_{\mathcal{I}}, p$} \Comment{function for new information}
\State $\vec{\textbf{s}}  \gets$ list of immediate successors of $a_t$ in $G$
\For{$\forall k \in$ $\vec{\textbf{s}}$}
\If {($rand(0, 1) \geq p$)}
\State node $k$ gains power
\State {$\vec{\textbf{y}}  \gets$}\Call{Simulation}{$k, G_{\mathcal{I}}, p$}\Comment{$\vec{y}$ stores the set of nodes that gains power}
\EndIf
\EndFor
\State \Return $\vec{\textbf{y}}$
\EndFunction

\end{algorithmic}
\end{minipage}}   

\end{algorithm}

\section{Flow Charts}
\label{App:FlowCharts}

\tikzstyle{startstop} = [rectangle, rounded corners, 
minimum width=1.5cm, 
minimum height=1cm,
text centered, 
draw=black, 
fill=white!30]

\tikzstyle{io} = [trapezium, 
trapezium stretches=true, 
trapezium left angle=70, 
trapezium right angle=110, 
minimum width=3cm, 
minimum height=1cm, text centered, 
draw=black, fill=white!30]

\tikzstyle{process} = [rectangle, 
minimum width=5cm, 
minimum height=1cm,  
text width=3cm, 
draw=black, 
fill=white!30]

\tikzstyle{decision} = [diamond, 
draw=black, 
fill=white!30]
\tikzstyle{arrow} = [thick,->,>=stealth]
\begin{figure}[h]
\centering

\begin{tikzpicture}[node distance=2cm, every node/.style={scale=0.59}]

\node (start) [startstop] {Start};

\node (pro1) [process, below of=start] {Pre-decision state $(S)$ : $(U^{+}, U^{-}_1, U^{-}_p, U^{-}_0)$};
\node (dec1) [decision, below of=pro1, yshift=-0.5cm] {Action $a_t$};
\node (pro2a) [process, below of=dec1, yshift=-0.5cm] {Post-decision state $(U^{a})$ : $(a, U^{+}, U^{-}_1, U^{-}_p, U^{-}_0)$};
\node (dec2) [decision, below of=pro2a, yshift=-1.0cm] {Aggregation? };
\node (pro2c) [process, below of=dec2, yshift=-1.5cm] {Aggregate Post-decision state using either $SA^1$, $SA^2$, or $SA^3$.
};
\node (in1) [io, right of=dec2, xshift=2cm] {$W$ revealed};
\node (pro2b) [process, right of=in1, xshift=2.5cm] {Pre-decision state $(\bar{S})$ : $(\bar{U^{+}}, \bar{U^{-}_1}, \bar{U^{-}_p}, \bar{U^{-}_0})$};
\node (dec3) [decision, below of=pro2b, yshift=-1.0cm] {$\bar{U^{-}_1} \cup \bar{U^{-}_p} = \phi$ };
\node (pro2d) [process, below of=dec3, yshift=-1.0cm] {Select an action from $\bar{U^{-}_1} \cup \bar{U^{-}_p}$};
\node (pro2e) [process, below of=pro2c, yshift=-1.0cm] {Store aggregated post-decision state in lookup table $\mathcal{L}$};
\node (stop) [startstop, left of=dec3, xshift= -1.25cm] {Stop};

\draw [arrow] (start) -- (pro1);

\draw [arrow] (pro1) -- (dec1);
\draw [arrow] (dec1) -- (pro2a);
\draw [arrow] (pro2a) -- (dec2);
\draw [arrow] (dec2) -- node[anchor=east] {Yes} (pro2c);
\draw [arrow] (dec2) -- node[anchor=south] {No} (in1);
\draw [arrow] (in1) -- (pro2b);
\draw [arrow] (pro2e) -| (in1);
\draw [arrow] (pro2b) -- (dec3);
\draw [arrow] (dec3) -- node[anchor=east] {Yes} (pro2d);
\draw [arrow] (dec3) -- node[anchor=south] {No} (stop);
\draw [arrow] (pro2c) -- (pro2e);
\draw [arrow] (pro2d.east) -- ++(1,0) |- (pro2a.east);
\draw [arrow] (dec3) -- (stop);

\end{tikzpicture}
\caption{Flowchart for solving a problem instance to completion to update approximate state values in the lookup table $\mathcal{L}$.} 
\label{fig:flowchart1}
\end{figure}
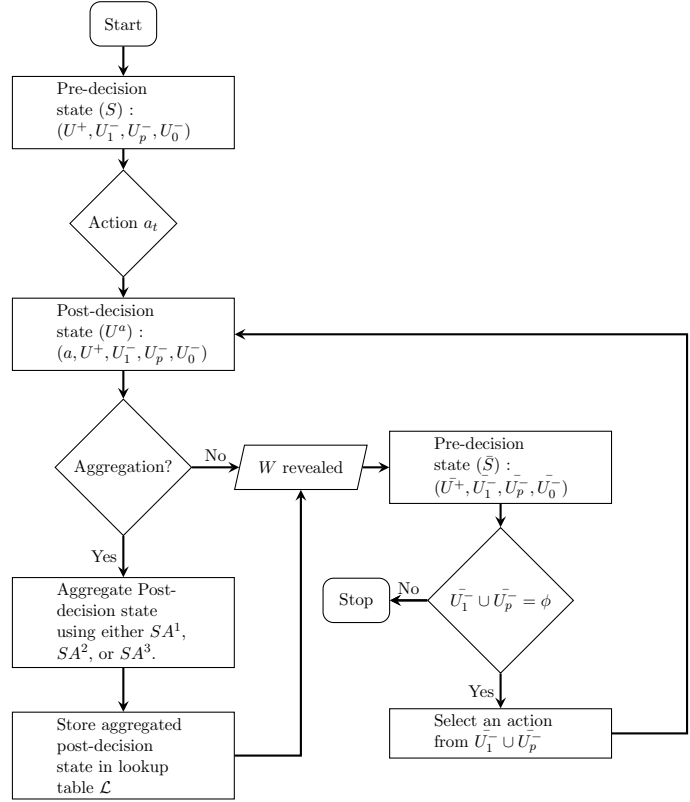

\begin{figure}
\centering
\begin{tikzpicture}[node distance=2cm, every node/.style={scale=0.59}]
\node (start) [startstop] {Start};
\node (in1) [io, right of=start, xshift=1.0cm] {Realization};
\node (pro1) [process, right of=in1, xshift=3.0cm] {Formulate the initial pre-decision state $(S)$};
\node (pro2a) [process, right of=pro1, xshift=4cm] {Select an action $a_t$ from $\mathcal{L}$ based on $S$ };
\node (pro2c) [process, right of=pro2a, xshift=4.5cm] {Update the post-decision state $S^{a}$};
\node (in2) [io, below of =pro2c, yshift=- 1.5cm] {Information $\omega$ revealed};
\node (pro2b) [process, left of=in2, xshift=-4.5cm] {Update the pre-decision state $\bar{S}$};
\node (dec3) [decision, left of=pro2b, xshift=-4.0cm] {$\bar{U^{-}_1} \cup \bar{U^{-}_p} = \phi$ };
\node (stop) [startstop, left of=dec3, xshift= -2.25cm] {Stop};
\draw [arrow] (start) -- (in1);
\draw [arrow] (in1) -- (pro1);
\draw [arrow] (pro1) -- (pro2a);
\draw [arrow] (pro2a) -- (pro2c);
\draw [arrow] (pro2c) -- (in2);
\draw [arrow] (in2) -- (pro2b);
\draw [arrow] (pro2b) -- (dec3);
\draw [arrow] (dec3.north)  -| node[anchor=west] {Yes} (pro2a.south);
\draw [arrow] (dec3) -- node[anchor=north] {No} (stop);
\end{tikzpicture}
\caption{Steps to compute the repair crew route for a given realization.}
\label{fig:flowchart2}
\end{figure}
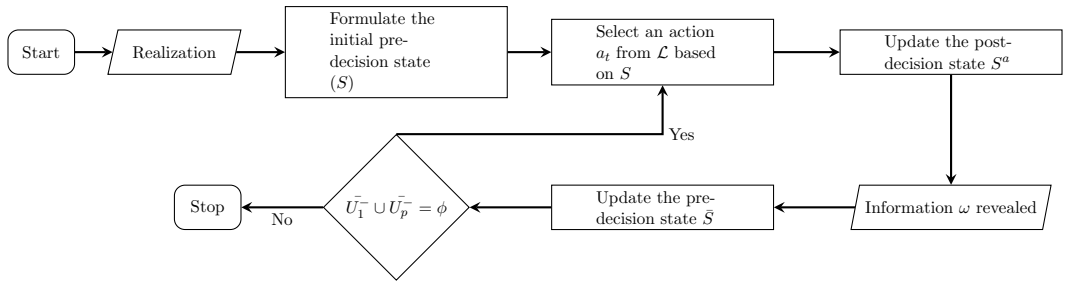

\section{Proof of Proposition \ref{Prop:PrereqRefinement}}
\label{Proof_of_Proposition_1}
\begin{proof}

To show the sub-optimality of solutions described in the proposition, we will show that we can find another solution with a lower objective function value under the conditions stated in the proposition.

Consider two routes for the TRNRP that are the same except for the sequence followed over nodes $j$ and $k$, where $k\preceq j$. We characterize these routes, called A and B, as follows.
\begin{center}
Route $A$: $\cdots i \rightarrow j \rightarrow k \rightarrow h \cdots$

Route $B$: $\cdots i \rightarrow k \rightarrow j \rightarrow h \cdots$
\end{center}
Let $C_A$ and $C_B$ represent the total expected service disruption time over all nodes following routes $A$ and $B$, respectively.
Note that these routes are determined offline; hence, during the crew travel, based on the information obtained from node status (faulty or not), the actual paths may differ. The list of all possible cases is shown in Table \ref{Tab:Cases} and explained below.
\begin{enumerate}[(i)]
\item After visiting node $i$, all nodes gain service.
\item After visiting node $i$, nodes $j$ and $k$ gain service, and there is at least one other node without service.
\item Node $k$ is not faulty, and all nodes gain service after visiting $j$.
\item Node $k$ is not faulty and gains service after visiting $j$, and there is at least one other node without service after visiting $j$.
\item Node $k$ is faulty and gains service after visiting $j$ and $k$, and all other nodes gain service after visiting $j$.
\item Node $k$ is faulty and does not gain service after visiting $j$, and there is at least one other node without service after visiting $j$.
\end{enumerate}

\begin{table}[H]
\begin{center}
 \scalefont{0.8}
\caption{ All possible actual paths based on routes $A$ and $B$.}
\begin{tabular}{c|c|c}
   \hline \hline
    Cases & Actual path following $A$ & Actual path following $B$  
    \\ \hline
    (i) & $\cdots i$  & $\cdots i$ \\
    \hline 
    (ii) & $\cdots i \rightarrow h \cdots$  & $\cdots i \rightarrow h \cdots$  \\
    \hline 
    (iii) & $\cdots i \rightarrow j \cdots$   & $\cdots i \rightarrow k \rightarrow j \cdots$  \\
    \hline 
    (iv) & $\cdots i \rightarrow j \rightarrow h \cdots$ & $\cdots i \rightarrow k \rightarrow j \rightarrow h \cdots$ \\
    \hline 
    (v) & $\cdots i \rightarrow j \rightarrow k \cdots$ &$\cdots i \rightarrow k \rightarrow j \cdots$ \\
    \hline 
    (vi) & $\cdots i \rightarrow j \rightarrow k \rightarrow h \cdots$   & $\cdots i \rightarrow k \rightarrow j \rightarrow h \cdots$ \\ 
    \hline \hline
\end{tabular}
    \label{Tab:Cases}
\end{center}
\end{table}
In cases (i) and (ii), the realized paths following routes $A$ and $B$ are the same. So, in such cases, $C_A = C_B$.

In cases (iii) and (iv), the path following solution $A$ avoids an unnecessary trip to node $k$, which occurs on the path following solution $B$. Hence, for those cases, $A$ is strictly better, i.e., $C_A < C_B$. 

In case (v), $Q_{i}^{\mathcal{A}} = Q_{i}^{\mathcal{B}}$ and $Q_{j}^{\mathcal{A}} \leq Q_{k}^{\mathcal{B}}$. Hence, $C_A\leq C_B$ if $d_{ij}\leq d_{ik}$.

In case (vi), $Q_{i}^{\mathcal{A}} = Q_{i}^{\mathcal{B}} = Q_{k}^{\mathcal{B}} \geq Q_{j}^{\mathcal{A}} \geq Q_{j}^{\mathcal{B}} =Q_{k}^{\mathcal{B}}$. That is, for each corresponding arc on the realized paths, the one on route $A$ maintains a smaller total number of nodes without service compared to
the one on route $B$. Hence, $C_A \leq C_B$, if $d_{ij} \leq d_{ik}$ and $d_{kh} \leq d_{jh}$.

Compiling all cases, solution $A$ cannot lead to a higher objective function value than solution $B$ if one of the following conditions holds:
\begin{itemize}
   \item $d_{ij} \leq d_{ik}$ and $d_{kh} \leq d_{jh}$ when $h \in V_{\mathcal{I}}$,
\item   $d_{ij} \leq d_{ik}$ when $h = 0$. 
\end{itemize}
Therefore, route $B$ is sub-optimal for the repair crew routing problem under these conditions.
\end{proof}

\section{Proof of Proposition 2}
 \label{Proof_of_Proposition_2}
\begin{proof}
    Consider an offline path $\mathcal{P}$.
    We know that $Q^{\mathcal{P}}_{1} \geq Q^{\mathcal{P}}_{2} \geq \dots \geq Q^{\mathcal{P}}_{|N_{I}|}$ since the number of nodes without service is non-increasing along the travel path.
    The corresponding expected total service disruption time is $\sum_{q=1}^{|N_I|} d_q^{\mathcal{P}} Q^{\mathcal{P}}_{q}$.
    
    Assume that path $\mathcal{P}$ is not double $k-opt$ guaranteed for some $k\geq 2$.    
    Let $\mathcal{P}^{'}$ be in the neighborhood of path $\mathcal{P}$ that can be obtained via a $k-opt$ local search move on path $\mathcal{P}$. Hence, there exist $q_1$ and $q_2$ such that $Q_{j}^{\mathcal{P}^{'}} < Q_{j}^{\mathcal{P}} $ for $j\geq  q_1$ and $Q_{j}^{\mathcal{P}^{'}} = Q_{j}^{\mathcal{P}} $ for $j<  q_1$; and $d_{k}^{\mathcal{P}^{'}} < d_{k}^{\mathcal{P}}$ and $d_{k}^{\mathcal{P}^{'}} = d_{k}^{\mathcal{P}}$ for $k \neq q_2$.
    Hence, the difference $\sum_{q=1}^{|N_I|} d_q^{\mathcal{P}} Q^{\mathcal{P}}_{q} - \sum_{q=1}^{|N_I|} d_q^{\mathcal{P}^{'}} Q^{\mathcal{P}^{'}}_{q}$ is positive, which implies that path $\mathcal{P}^{'}$ leads to a better objective function value.
\end{proof}

\section{Instance Generation}
\label{Ap:InstanceGeneration}

In this section, we present our instance generation in detail. We create our test instances in two steps.
In the first step, we determine $N_{\mathcal{I}}$ random demand node locations on a given region with different shapes, e.g., rectangular or circular. The location of the depot (node 0) is chosen as the center of gravity. The road network includes edges between every pair of nodes, and the travel times are equal to the edge lengths. 
We present the pseudo-code to generate the random node locations in Algorithm \ref{alg:location1}. 

After generating the locations of the nodes on the service region, we create the arcs on the infrastructure network ($G_{\mathcal{I}}$) by forming a minimum spanning tree with $N_{\mathcal{I}}-1$ arcs. 
We later modify the initially created networks to obtain instances with specific attributes using a degree-constrained minimum spanning tree ($d$-MST) heuristic algorithm, where the degree of a node cannot exceed a certain value, which is referred to as the $degree$ constraint.
The heuristic algorithm we use is based on the algorithm proposed in \cite{boldon1996minimum}.
Let $l_{max}$ and $l_{min}$ be the maximum and the minimum edge lengths in the original MST. 
In the $d$-MST heuristic algorithm, we first find any nodes violating the degree constraint, i.e., nodes whose degree exceeds the value $degree$. Let the set of such nodes be $F$. Then, for each node $i\in F$, we check if node $i$ has an adjacent edge $e_i$ on the MST that is longer than $l_{min}$, i.e., $l_{e_i} > l_{min}$. If so, we update the length of edge $e_i$ as follows
\begin{equation*}
   l_{e_{i}}' = l_{e_{i}} + f\cdot l_{max}\frac{l_{e_{i}} - l_{min}}{l_{max} - l_{min}},
\end{equation*}
where $f$ is the total number of nodes in the current MST violating the degree constraint. We perform this for each edge $e_i$ such that $l_{e_i} > l_{min}$. By doing so, we penalize the edges that are incident to nodes with degree violations. 
We then recalculate the MST using the updated edge lengths until all nodes satisfy the degree bound. The pseudo-code of this heuristic is presented in Algorithm \ref{algo:dcmst}.
Once the degree-constrained MST is obtained, we determine the node indices following the procedure in Algorithm \ref{alg:rename}. We designate the node with the highest degree as the source node (node 1).

To generate power network topologies with varying depths, we randomly select $\frac{N_{\mathcal{I}}}{2}$ and $\frac{N_{\mathcal{I}}}{4}$ nodes from the obtained $d$-MST and connect them to the grandparent nodes avoiding any intersection of the edges. 
\begin{algorithm}[H]
    \caption{ Generating Node Coordinates for Different Shapes of Service Regions}\label{alg:location1}
\resizebox{0.8\textwidth}{!}{
\begin{minipage}{\textwidth}    
    \begin{algorithmic}[1]

\Function{Place\_Points\_Inside\_Rectangle}{$w, h, N_{\mathcal{I}}$}
    
    \For{$i$ \textbf{in} $1$ to $N_{\mathcal{I}}$}
        \State $x \gets \text{rn.uniform}(X_0, X_0 + w)$
        \State $y \gets \text{rn.uniform}(Y_0, Y_0 + h)$
        \State add $x,y$ to $X, Y$, respectively.
    \EndFor
    \State \textbf{return} $X, Y$
\EndFunction 

\Function{Place\_Points\_Inside\_Circle}{$R, N_{\mathcal{I}}$}
    
        \For{$i$ \textbf{in} $\text{range}(N_{\mathcal{I}})$}
            \State $r \gets \text{rn.uniform}(0, R)$
            \State$\alpha \gets \text{rn.uniform}(0, 2\pi)$
            \State $x \gets X_0 + r \times \cos(\alpha)$
            \State $y \gets Y_0 + r \times \sin(\alpha)$
            \State add $x,y$ to $X, Y$, respectively.
    \EndFor
    \State \textbf{return} $X, Y$
\EndFunction 

\If{($flag_{location} == 0$)}
               \State $X, Y \gets \Call{Place\_Points\_Inside\_Rectangle}{w, h, N_{\mathcal{I}}}$
            \Else
                 \State $X, Y \gets
                \Call{Place\_Points\_Inside\_Circle}{R, N_{\mathcal{I}}}$

        \EndIf

\end{algorithmic}
\end{minipage}}
\end{algorithm}

\begin{algorithm}[H]
\caption{$d$-MST Heuristic}
\label{algo:dcmst}
\resizebox{0.8\textwidth}{!}{
\begin{minipage}{\textwidth}
\begin{algorithmic}[1]
\Function{is\_degree\_satisfied}{$mst, degree$}
    \For{$i$ \textbf{in all the nodes of} $mst$ }
        \If{$\text{degree of node } i > degree$}
            \State \Return \textbf{False}
        \EndIf
    \EndFor
    \State \Return \textbf{True}
\EndFunction

\Function{get\_edges\_adjacent\_to\_violating\_node}{$mst, node$}
    \State $adjacent\_edges \gets []$
    \For{$u, v$ in $mst$}
        \If{$u = node$ \textbf{or} $v = node$}
            \State add the edge $u,v$ along with $d_{uv}$ in $adjacent\_edges$
        \EndIf
    \EndFor
    \State \Return $adjacent\_edges$
\EndFunction

\Function{find\_lengths\_greater\_than\_least\_edge\_in\_mst}{$adjacent\_edges$}
    \State $sorted\_edges$ $\gets$ store the edges whose length is greater than the minimum length
    \State \Return $sorted\_edges$
\EndFunction

\Function{degree\_constrained\_mst\_heuristics}{$G_{\mathcal{I}}, degree$}
    \State $current\_graph \gets \text{create a copy of the } G_{\mathcal{I}}$
    \State $mst \gets \text{find the minimum spanning tree of the } current\_graph$
    \While{\textbf{not} \Call{is\_degree\_satisfied}{\textit{mst}, \textit{degree}}}
        \State $violating\_nodes \gets \text{stores the nodes that violates the degree condition in } mst$
        \State $l\_min \gets \text{obtains the minimum length in } mst$
        \State $l\_max \gets \text{obtains the maximum length in } mst$
        \For{$k$ \textbf{in} $violating\_nodes$}
            \State $adjacent\_edges\gets \Call{get\_edges\_adjacent\_to\_violating\_node}{mst, k}$
            \State $lengths \gets \Call{find\_lengths\_greater\_than\_least\_edge\_in\_mst}{adjacent\_edges}$
            \For{$u, v, l$ \textbf{in} $adjacent\_edges$}
                \If{$l \in lengths$}
                    \State add penalty to the edges and update the length of the edges
                \EndIf
            \EndFor
        \EndFor
        \State $mst \gets$ update the minimum spanning tree of the $current\_graph$
    \EndWhile
    \State \Return $mst$
\EndFunction
\State $mst_{edges} \gets \text{\Call{degree\_constrained\_mst\_heuristics}{$G_{\mathcal{I}}, degree$}}$
\end{algorithmic}
\end{minipage}}

\end{algorithm}


\begin{algorithm}[H]
\caption{Renaming Nodes in the Graph Obtained by $d$-MST}\label{alg:rename}
\resizebox{0.8\textwidth}{!}{
\begin{minipage}{\textwidth}
\begin{algorithmic}[1]
\State create the $adjacency\_list$ of the MST obtained from $d$-MST
\State $new\_name \gets$ select randomly a node and assign it $1$
\State $connected\_nodes \gets$ stores list of nodes having an edge with the recently modified node number
\State $nodes\_renamed \gets$ stores the modified node numbers
\Function{rename}{$connected\_nodes, new\_name$}
    \If{len($connected\_nodes$) != 0}
        \For{$node$ in $connected\_nodes$}
            \State $h \gets node$
            \State  $name \gets new\_name + 1$ 
            \State add $name$, and update $node$ with $name$ everywhere in the $adjacency\_list$.
            \State $new\_name \gets name$ and add it to $nodes\_renamed$
        \EndFor
        \State update $connected\_nodes$, $new\_name$ and decrease the length of $nodes\_renamed$
        \If{len($new\_name$) = $N$}
            \State \Return $adjacency\_list$
        \EndIf
        \State \Call{rename}{$connected\_nodes, new\_name$}
    \Else
        \State update $connected\_nodes$, $new\_name$ and decrease the length of $nodes\_renamed$
        \If{len($new\_name$) = $N$}
            \State \Return $adjacency\_list$
        \EndIf
        \State \Call{rename}{$connected\_nodes, new\_name$}
    \EndIf
    \State \Return \textit{adjacency\_list}
\EndFunction
\State s = \Call{rename}{\textit{connected\_nodes}, \textit{new\_name}}
\State convert the $adjacency\_list$ back to the graph.
\end{algorithmic}
\end{minipage}}
\end{algorithm}

\section{Extended Computational Results}
\label{App:ExtendedResults}

\begin{figure}[H]
\centering
\subfloat[Original network with depth $d=15$.]{\label{fig40_1}{\includegraphics[width=0.4\textwidth]{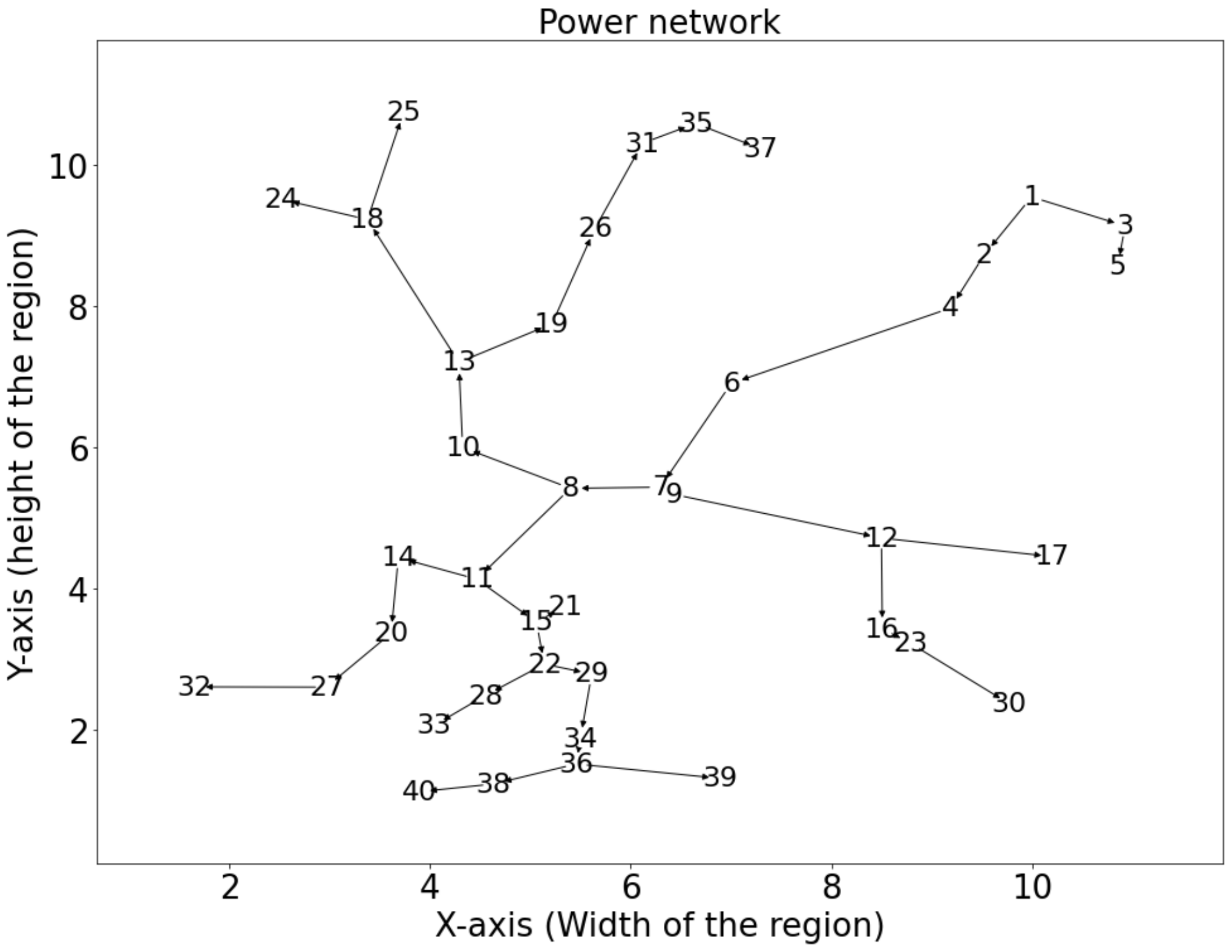}}}\hspace{1cm}
\subfloat[Modified network with depth $d=9$.]{\label{fig40_2}{\includegraphics[width=0.40\textwidth]{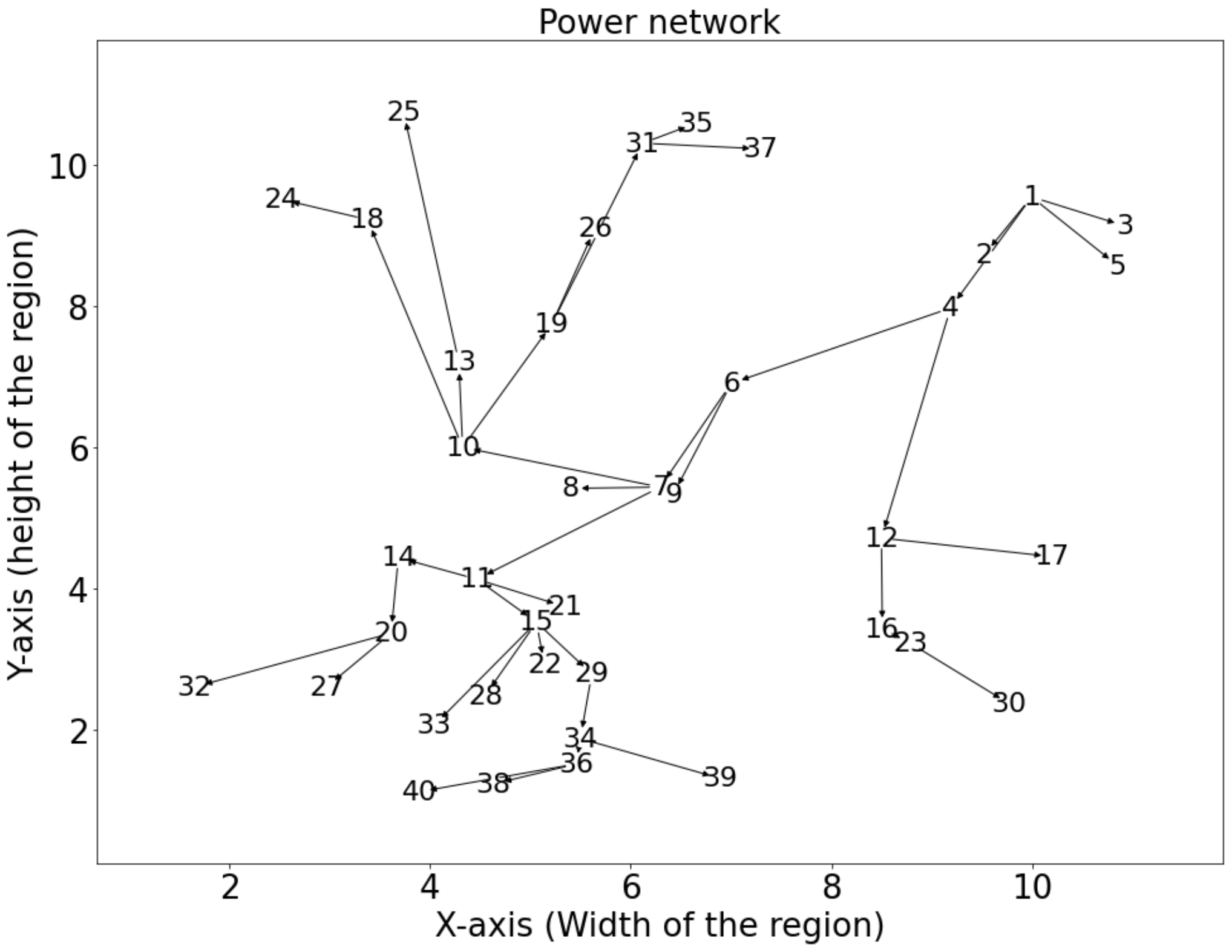}}}
\caption{Node dispersions of power distribution networks with 40 nodes on a 10$\times$10 square region.}
\label{Fig:40_Square}
\end{figure}

\begin{table}[H]
\centering
\caption{Comparison of average (Avg) objective values and \% CI's for gap compared to $SA^3$ solutions for 40 nodes on a 10$\times$10 square region for different depths.}
\label{Tab:40_sq}
\resizebox{0.8\textwidth}{!}{
\begin{tabular}{l|l|r|rr|r|rr}
\hline \hline 
\multicolumn{2}{c|}{} & \multicolumn{3}{c|}{\textbf{Network depth} $d=15$} 
& \multicolumn{3}{c}{\textbf{Network depth} $d=9$} \\ \hline
\multicolumn{2}{l|}{} & \multicolumn{3}{c|}{\makecell{\textbf{Avg Total Service } \\ \textbf{Disruption Times}}} 
& \multicolumn{3}{c}{\makecell{\textbf{Avg Total Service} \\ \textbf{ Disruption Times}}}\\ \hline
      $s$ & $p$ & \multicolumn{1}{c|}{\textit{Our Methods}} &  \multicolumn{2}{c|}{\textit{Benchmarks}}   &
     \multicolumn{1}{c|}{\textit{Our Methods}}  &  \multicolumn{2}{c}{\textit{Benchmarks}}\\ &  & $SA^3$ & $PS$ & $NN$ 
& $SA^3$ & $PS$ & $NN$ \\ \hline

\multirow{4}{*}{0}   & 0.25 & 533  & 646  & 1379  
& 462  & 580  & 1208 \\
                     & 0.5  & 809  & 1296 & 1599 
                     & 759  & 1308 & 1495  \\
                     & 0.75 & 1006 & 2142 & 1692  
                     & 992  & 2294 & 1669 \\
                     & 0.9  & 1087 & 2718 & 1725 
                     & 1076 & 2978 & 1736  \\ \hline \multicolumn{2}{l|}{\multirow{1}{*}{\textbf{CI for Avg \% Gap }}} & -  & (87, 89)  & (98, 101)  
                     &
                    -   & (102, 104)  & (106, 108) \\ \hline
                     
\multirow{4}{*}{1.5} & 0.25 & 823  & 922  & 1880  
& 704  & 796  & 1665   \\
                     & 0.5  & 1406 & 1866 & 2599  
                     & 1314 & 1792 & 2453  \\
                     & 0.75 & 1926 & 3033 & 3166  
                     & 1883 & 3119 & 3129 \\
                     & 0.9  & 2198 & 3811 & 3473  
                     & 2180 & 4039 & 3491   \\ \hline
                      \multicolumn{2}{l|}{\multirow{1}{*}{\textbf{CI for Avg \% Gap }}} & -  & (44, 45)  & (86, 88)    
                      & -   & (49, 50)  & (94, 96) \\ \hline
\multirow{4}{*}{3}   & 0.25 & 1110 & 1199 & 2381  
& 932  & 1013 & 2122 \\
                     & 0.5  & 2000 & 2436 & 3598  
                     & 1854 & 2276 & 3410 \\
                     & 0.75 & 2838 & 3925 & 4639 
                     & 2784 & 3944 & 4590  \\
                     & 0.9  & 3307 & 4904 & 5221 
                     & 3278 & 5101 & 5245  \\ \hline
                      \multicolumn{2}{l|}{\multirow{1}{*}{\textbf{CI for Avg \% Gap }}} & -  & (30, 31)  & (80, 82)   
                      & -   & (32, 33)  & (89, 91) \\ \hline
 \hline 
\end{tabular}} 
\end{table}

\begin{figure}[H]
\centering
\subfloat[Original network with depth $d=22$.]{\label{fig40_1}{\includegraphics[width=0.40\textwidth]{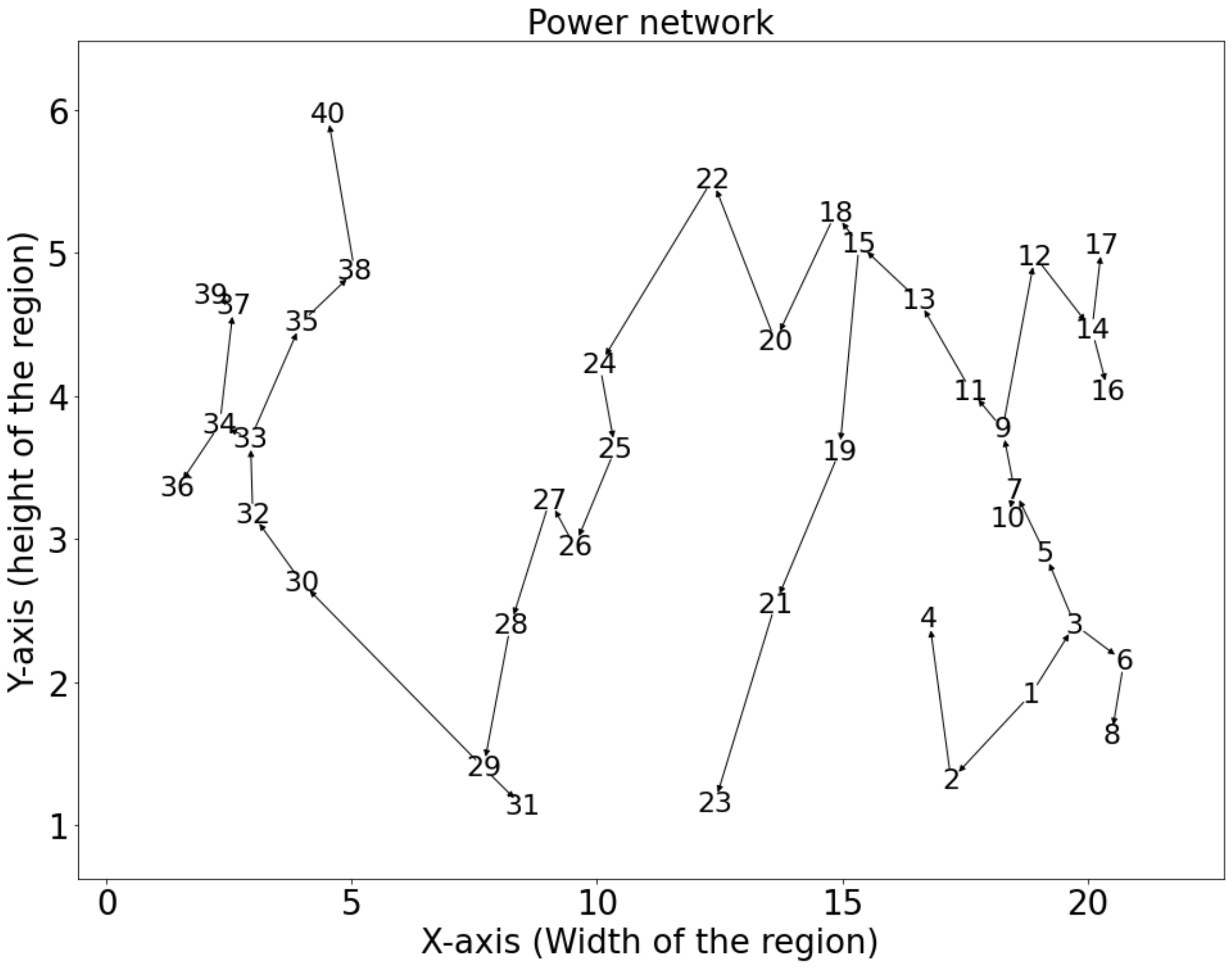}}}\hspace{1cm}
\subfloat[Modified network with depth $d=12$.]{\label{fig40_2}{\includegraphics[width=0.40\textwidth]{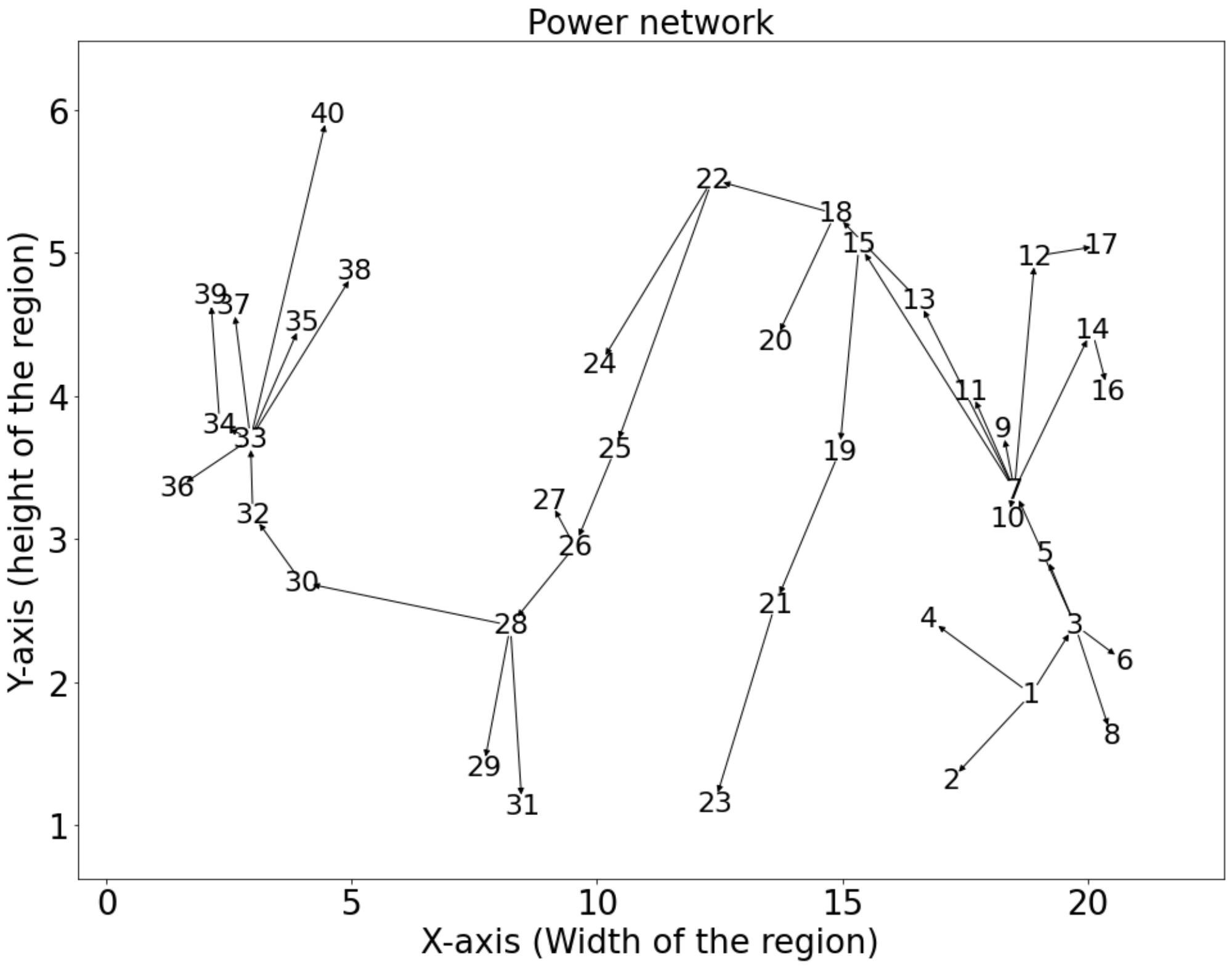}}}
\caption{Node dispersions of power distribution networks with 40 nodes on a 20 $\times$ 5 rectangle region.}
\label{Fig:40_rec}
\end{figure}

\begin{table}[H] 
\centering 
\caption{Comparison of average (Avg) objective values and \% CI's for gap compared to $SA^3$ solutions for 40 nodes on a 20 $\times$ 5 rectangle region for different depths.} 
\label{Tab:40_rec}

\resizebox{0.8\textwidth}{!}{
\begin{tabular}{l|l|r|rr|r|rr}
\hline \hline 
\multicolumn{2}{c|}{} & \multicolumn{3}{c|}{\textbf{Network depth} $d=22$} 
& \multicolumn{3}{c}{\textbf{Network depth} $d=12$} \\ \hline
\multicolumn{2}{l|}{} & \multicolumn{3}{c|}{\makecell{\textbf{Avg Total Service } \\ \textbf{Disruption Times}}} 
& \multicolumn{3}{c}{\makecell{\textbf{Avg Total Service} \\ \textbf{ Disruption Times}}}\\ \hline
      $s$ & $p$ & \multicolumn{1}{c|}{\textit{Our Methods}} &  \multicolumn{2}{c|}{\textit{Benchmarks}}   & 
     \multicolumn{1}{c|}{\textit{Our Methods}}  &  \multicolumn{2}{c}{\textit{Benchmarks}}\\ &  & $SA^3$ & $PS$ & $NN$ 
& $SA^3$ & $PS$ & $NN$ \\ \hline
\multicolumn{1}{l|}{\multirow{4}{*}{0}} & 0.25 & 935 & 1185 & 2290  
& 874 & 1115 & 2280 \\ 
\multicolumn{1}{l|}{} & 0.5 & 1156 & 2065 & 2369  
& 1128 & 2112 & 2362 \\ 
\multicolumn{1}{l|}{} & 0.75 & 1285 & 3114 & 2423 
& 1213 & 3381 & 2419 \\ 
\multicolumn{1}{l|}{} & 0.9 & 1362 & 3941 & 2447 
& 1313 & 4359 & 2445  \\ 
\hline 
  \multicolumn{2}{l|}{\multirow{1}{*}{\textbf{CI for Avg \% Gap }}} & -  & (111, 114)  & (109, 108)  
  & -   & (124, 126)  & (114, 116) \\ \hline
\multicolumn{1}{l|}{\multirow{4}{*}{1.5}} & 0.25 & 1252 & 1461 & 2906  
& 1137 & 1332 & 2887 \\ 
\multicolumn{1}{l|}{} & 0.5 & 1778 & 2627 & 3565  
& 1733 & 2588 & 3551 \\ 
\multicolumn{1}{l|}{} & 0.75 & 2228 & 3994 & 4194  
& 2207 & 4194 & 4186\\ 
\multicolumn{1}{l|}{} & 0.9 & 2488 & 5030 & 4564  
&  2458 & 5414 & 4560 \\ 
\hline 
  \multicolumn{2}{l|}{\multirow{1}{*}{\textbf{CI for Avg \% Gap }}} & -  & (62, 64)  & (104, 106) 
  & -   & (69, 71)  & (110, 111) \\ \hline
\multicolumn{1}{l|}{\multirow{4}{*}{3}} & 0.25 & 1560 & 1736 & 3522 
& 1383 & 1548 & 3494 \\ 
\multicolumn{1}{l|}{} & 0.5 & 2389 & 3189 & 4761  
& 2297 & 3064 & 4740 \\ 
\multicolumn{1}{l|}{} & 0.75 & 3155 & 4873 & 5966 
& 3132 & 5006 & 5952 \\ 
\multicolumn{1}{l|}{} & 0.9 & 3598 & 6118 & 6681  
& 3572 & 6470 & 6676 \\ 
\hline 
  \multicolumn{2}{l|}{\multirow{1}{*}{\textbf{CI for Avg \% Gap }}} & -  & (41, 43)  & (101, 102)  
  & -   & (46, 48)  & (111, 112) \\ 
\hline \hline 
\end{tabular}}
\end{table}

\begin{figure}[H]
\centering
\subfloat[Original network with depth $d=21$.]{\label{fig40_1}{\includegraphics[width=0.40\textwidth]{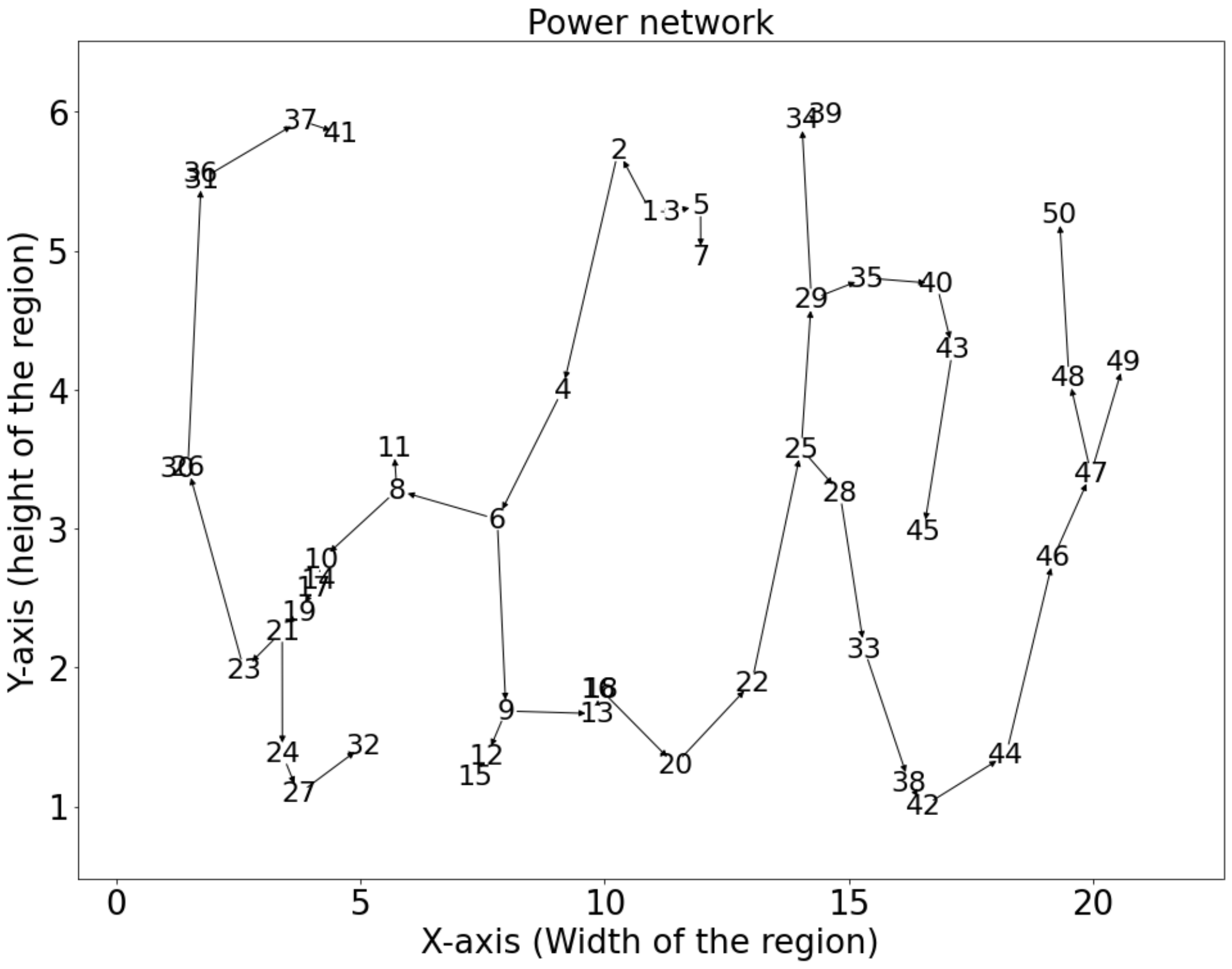}}}\hspace{1cm}
\subfloat[Modified network with depth $d=10$.]{\label{fig40_2}{\includegraphics[width=0.40\textwidth]{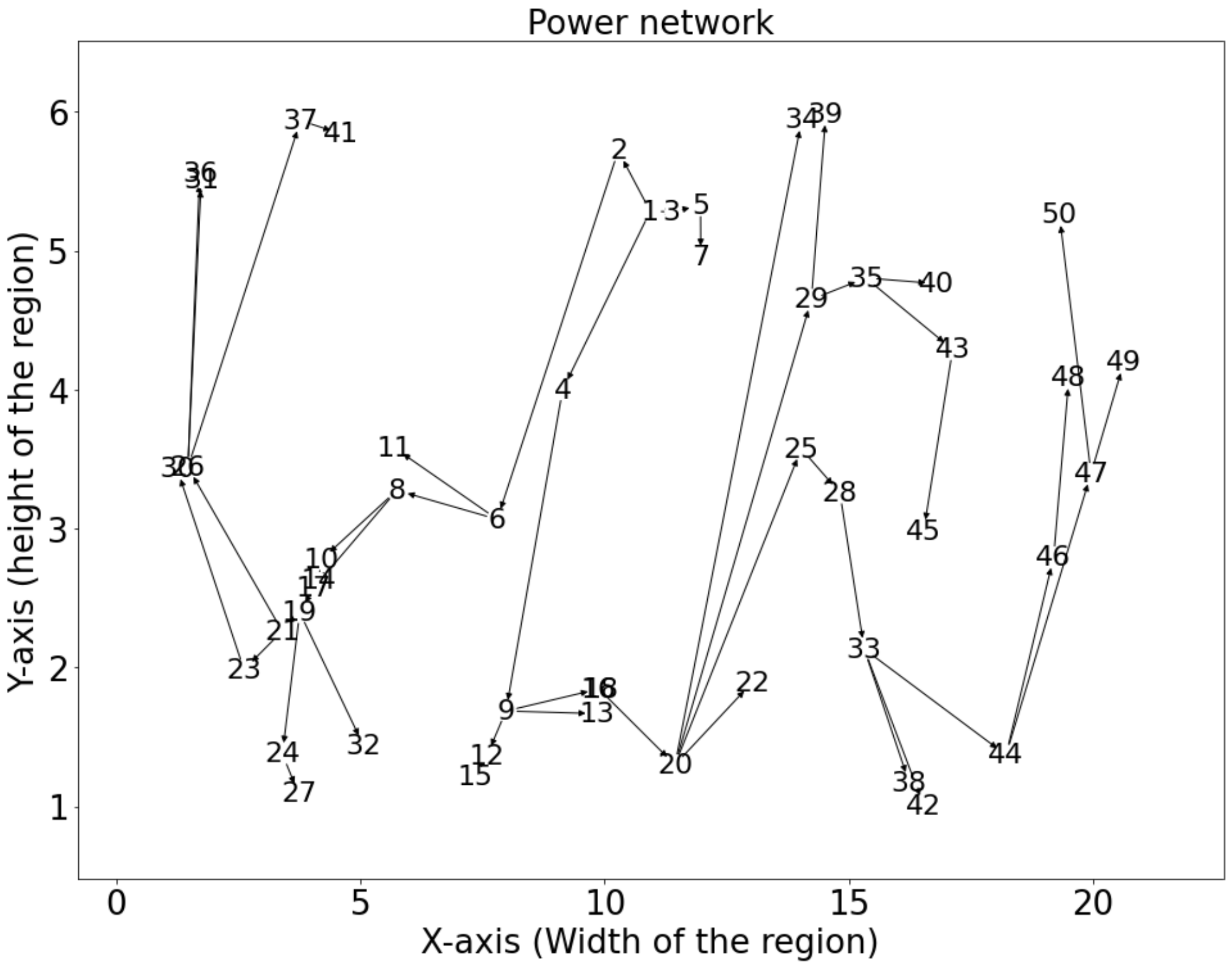}}}
\caption{Node dispersions of power distribution networks with 50 nodes on a 20 $\times$ 5 rectangle region.}
\label{Fig:50_rec}
\end{figure}

\begin{table}[H]

\centering
\caption{Comparison of average (Avg) objective values and \% CI's for gap compared to $SA^3$ solutions for 50 nodes on a 20 $\times$ 5 rectangle region for different depths.}
\label{Tab:50_rec}

\resizebox{0.8\textwidth}{!}{
\begin{tabular}{l|l|r|rr|r|rr}
\hline \hline 
\multicolumn{2}{c|}{} & \multicolumn{3}{c|}{\textbf{Network depth} $d=21$}
& \multicolumn{3}{c}{\textbf{Network depth} $d=10$} \\ \hline
\multicolumn{2}{l|}{} & \multicolumn{3}{c|}{\makecell{\textbf{Avg Total Service } \\ \textbf{Disruption Times}}} 
& \multicolumn{3}{c}{\makecell{\textbf{Avg Total Service} \\ \textbf{ Disruption Times}}}\\ \hline
      $s$ & $p$ & \multicolumn{1}{c|}{\textit{Our Methods}} &  \multicolumn{2}{c|}{\textit{Benchmarks}}   & 
     \multicolumn{1}{c|}{\textit{Our Methods}}  &  \multicolumn{2}{c}{\textit{Benchmarks}}\\ &  & $SA^3$ & $PS$ & $NN$ 
& $SA^3$ & $PS$ & $NN$ \\ \hline
\multicolumn{1}{l|}{\multirow{4}{*}{0}}   & 0.25     & 1126 & 1745  & 3501  
& 1019 & 1518  & 3464  \\
\multicolumn{1}{l|}{}                     & 0.5      & 1489 & 3765  & 3744    
& 1417 & 3540  & 3725   \\
\multicolumn{1}{l|}{}                     & 0.75     & 1714 & 6764  & 3906   
& 1708 & 6525  & 3898   \\
\multicolumn{1}{l|}{}                     & 0.9      & 1840 & 8427  & 3981   
& 1784 & 9319  & 3978   \\ \hline
 \multicolumn{2}{l|}{\multirow{1}{*}{\textbf{CI for Avg \% Gap }}} & -  & (211, 214)  &  (154, 156)  
 & -   & (192, 193)  & (164, 166) \\ \hline
\multicolumn{1}{l|}{\multirow{4}{*}{1.5}} & 0.25     & 1587 & 2137  & 4428    
& 1415 & 1821  & 4367  \\
\multicolumn{1}{l|}{}                     & 0.5      & 2404 & 4597  & 5539    
& 2391 & 4241  & 5504   \\
\multicolumn{1}{l|}{}                     & 0.75     & 3145 & 7869  & 6589    
& 3121 & 8007  & 6573  \\ 
\multicolumn{1}{l|}{}                     & 0.9      & 3591 & 10107 & 7195     
& 3526 & 10948 & 7188     \\ \hline
 \multicolumn{2}{l|}{\multirow{1}{*}{\textbf{CI for Avg \% Gap }}} & -  & (114, 116)  & (132, 132)  
 & -   & (105, 106)  & (138, 140) \\ \hline
\multicolumn{1}{l|}{\multirow{4}{*}{3}}   & 0.25     & 2046 & 2528  & 5355   
& 1752 & 2124  & 5270   \\
\multicolumn{1}{l|}{}                     & 0.5      & 3325 & 5428  & 7334     
& 3221 & 4943  & 7284 \\
\multicolumn{1}{l|}{}                     & 0.75     & 4775 & 9214  & 9272     
& 4575 & 9250  & 9247   \\
\multicolumn{1}{l|}{}                     & 0.9      & 5546 & 11786 & 10408    
& 5245 & 12576 & 10399  \\ \hline
\multicolumn{2}{l|}{\multirow{1}{*}{\textbf{CI for Avg \% Gap }}} & -  & (78, 79)  & (123, 126) 
& -   & (68, 68)  & (129, 130) \\ \hline
\hline
\end{tabular}} 
\end{table}

\end{document}